\documentclass[DIV=16,twocolumn,letterpaper,headings=normal]{scrartcl}       % twocolumn (second format)
\usepackage{mathptmx}
\usepackage{amsmath,amssymb,graphics,subfig,url,tikz}
\usepackage[utf8]{inputenc}
\usepackage[english]{babel}
\usepackage[pdftex,bookmarks]{hyperref}

\hypersetup{pdftitle = {A Study on Using Hierarchical Basis Error Estimates in Anisotropic
      Mesh Adaptation for the Finite Element Method},
   pdfauthor = {Lennard Kamenski},
   pdfsubject = {MSC 2010: 65N30, 65N50},
   pdfkeywords = {mesh adaptation, anisotropic mesh, finite element,
   a~posteriori estimate, hierarchical basis, variational problem,
   anisotropic diffusion}
}

\providecommand{\boldx}{\boldsymbol{x}}
\providecommand{\R}{\mathbb{R}}
\providecommand{\cT}{\mathcal{T}}
\providecommand{\abs}[1]{\left|{#1}\right|}
\providecommand{\norm}[1]{\left\|{#1}\right\|}

\begin{document}

%*** Front Matter **************************************************

%--- title ---------------------------------------------------------
\title{A Study on Using Hierarchical Basis Error Estimates in Anisotropic
      Mesh Adaptation\\ for the Finite Element Method}

\author{Lennard Kamenski
   \thanks{Department of Mathematics, University of Kansas, $\qquad\qquad$
       Email: \nolinkurl{lkamenski@math.ku.edu}.} ~
   \thanks{This research was supported in part by the German 
      Research Foundation through the grant KA~3215/1-1.}}

\date{\normalsize{June 28, 2011}}

%\titlerunning{Using Hierarchical Error Estimates 
%   for Anisotropic Mesh Adaptation}

\maketitle
%--- abstract -------------------------------------------------------
\noindent \textbf{Abstract:} 
A common approach for generating an anisotropic mesh is the $M$-\emph{uniform mesh} approach where an adaptive mesh is generated as a uniform one in the metric specified by a given tensor $M$.
A key component is the determination of an appropriate metric which is often based on some type of Hessian recovery.
Recently, the use of a global hierarchical basis error estimator was proposed for the development of an anisotropic metric tensor for the adaptive finite element solution.
This study discusses the use of this method for a selection of different applications.
Numerical results show that the method performs well and is comparable with existing metric tensors based on Hessian recovery.
Also, it can provide even better adaptation to the solution if applied to problems with gradient jumps and steep boundary layers.
For the Poisson problem in a domain with a corner singularity, the new method provides meshes that are fully comparable to the theoretically optimal meshes.



%*** Paper Content **************************************************

\section{Introduction}
\label{sec:introduction}
%********************************************************************
A common approach for generating an anisotropic mesh is the $M$-\emph{uniform mesh} approach based on generation of a quasi-uniform mesh in the metric space defined by a symmetric and strictly positive definite metric tensor $M$. 
A scalar metric tensor will lead to an isotropic mesh while a full metric tensor will generally result in an anisotropic mesh.
In this sense, the mesh generation procedure is the same for both isotropic and anisotropic mesh generation.
A key component of the approach is the determination of an appropriate metric often based on some type of error estimates.

Typically, the appropriate metric tensor depends on the Hessian of the exact solution of the underlying problem, which is often unavailable in practical computation, thus requiring the recovery of an approximate Hessian from the computed solution.
A number of recovery techniques are used for this purpose, for example the gradient recovery technique by Zienkiewicz and Zhu \cite{ZieZhu92,ZieZhu92a}, the technique based on the variational formulation by Dolej\v{s}{\'\i} \cite{Dolejs98}, or the quadratic least squares fitting (QLS) proposed by Zhang and Naga \cite{ZhaNag05}.
Generally speaking, Hessian recovery methods work well when exact nodal function values are provided (e.g. interpolation problems), but unfortunately they do not provide an accurate recovery when applied to linear finite element approximations on non-uniform meshes, as pointed out by the author in \cite{Kamens09PhD}.
Recently, conditions for asymptotically exact gradient and convergent Hessian recovery from a hierarchical basis error estimator have been given by Ovall \cite{Ovall07}. 
His result is based on superconvergence results by Bank and Xu \cite{BanXu03,BanXu03a}, which  require the mesh to be uniform or almost uniform: assumptions which are usually violated by adaptive meshes.

Hence, a convergence of adaptive algorithms based explicitly on the Hessian recovery cannot be proved in a direct way, even if their application is successful in practical computations \cite{Dolejs98,HuaLi10,LiHua10,VasLip99}.
This explains the recent interest in anisotropic adaptation strategies  based on some type of \emph{a~posteriori} error estimates.
For example, Cao et al. \cite{CaHuRu01} studied two \emph{a~posteriori} error estimation strategies for computing scalar monitor functions for use in adaptive mesh movement; Apel et al. \cite{ApGrJM04} investigated a number of \emph{a~posteriori} strategies for computing error gradients used for directional refinement; and Agouzal et al. \cite{AgLiVa08,AgLiVa09,AgLiVa10} and Agouzal and Vassilevski \cite{AgoVas10} proposed a new method for computing metric tensors to minimize the interpolation error provided that an edge-based error estimate is given.

Recently, Huang et al. \cite{HuKaLa10} presented a mesh adaptation method based on hierarchical basis error estimates (HBEE).
The new framework is developed for the linear finite element solution of a boundary value problem of a second-order elliptic partial differential equation (PDE), but it is quite general and can easily be adopted to other problems. 
A key idea in the new approach is the use of the globally defined HBEE for the reliable directional information: globally defined error estimators have the advantage that they contain more directional information of the solution; error estimation based on solving local error problems, despite its success in isotropic mesh adaptation, do not contain enough directional information, which is global in nature; moreover, Dobrowolski et al. \cite{DoGrPf99} have pointed out that local error estimates can be inaccurate on anisotropic meshes.
A brief description of the method is provided in Sect.~\ref{sec:metricBasedAdaptation}.

The objective of this article is to study the application of this new anisotropic adaptation approach to different problems.

The first example deals with a boundary value problem of a second-order elliptic PDE.
Numerical results in \cite{HuKaLa10} have already shown that the new mesh adaptation approach can be a successful alternative to Hessian recovery in mesh adaptation for a boundary value problem of a second-order elliptic PDE.
In this paper, we would like to investigate another interesting question, namely whether the proposed algorithm is able to provide optimal mesh adaptation for a problem where the optimal mesh design is known theoretically.
For this purpose, a Poisson problem in a domain with a corner singularity is studied, which requires proper mesh adaptation in order to obtain an accurate finite element solution.
Theoretically optimal meshes are known for this problem and are based on some \emph{a~priori} information of the domain and the solution.
Section~\ref{sec:cornerSingularity} investigates if the proposed method is able to generate adaptive meshes that are comparable to the optimal ones without any \emph{a~priori} information on the exact solution.
A brief overview over the considered mathematical problem, the properties of the optimal meshes and the numerical results for the new method are provided.

Section~\ref{sec:variationalProblem} presents an anisotropic metric tensor for general variational problems developed by Huang et al. \cite{HuKaLi10} using the HBEE and the underlying variational formulation and gives a numerical example for a non-quadratic variational problem.
The metric tensor is completely \emph{a~posteriori}: it is based solely on the residual, edge jumps, and the \emph{a~posteriori} error estimate.

The third example is an anisotropic diffusion problem. The exact solution of this problem satisfies the maximum principle and it is desirable for the numerical solution to fulfill its discrete counterpart: the discrete maximum principle (DMP).
Recently, Li and Huang \cite{LiHua10} developed an anisotropic metric tensor based on the anisotropic non-obtuse angle condition, which provides both mesh adaptation and DMP satisfaction for the numerical solution: the mesh alignment is determined by the main diffusion drag direction, i.e. by the underlying PDE, and the Hessian of the exact solution determines the optimal mesh density.
In Sect.~\ref{sec:anisotopicDiffusion}, the Hessian of the exact solution is replaced with the Hessian of the hierarchical error estimator to obtain a new, completely \emph{a~posteriori}, metric tensor accounting for both DMP satisfaction and mesh adaptation.

Concluding remarks on the numerical examples and some key components of the hierarchical basis error estimator are given in Sect.~\ref{sec:concludingRemarks}.

\section{Anisotropic mesh adaptation based on hierarchical basis error estimator}
\label{sec:metricBasedAdaptation}
%********************************************************************
Consider the solution of a variational problem: find $u \in V$ such that
\begin{align}
   a(u,v) = f(v) \quad \forall v \in V
   \tag{$P$}
   \label{eq:vp_H}
\end{align}
where $V$ is an appropriate Hilbert space of functions over a domain $\Omega \in \R^2$, $a(\cdot,\cdot)$ is a bilinear form defined on $V\times V$, and $F(\cdot)$ is a continuous linear functional on $V$.

For a given simplicial mesh $\cT_h$, the linear finite element approximation $u_h$ of $u$ is the solution of the corresponding variational problem in a finite dimensional subspace $V_h \subset V$ of piecewise linear functions over $\cT_h$: find  $u_h\in V_h$ such that
\begin{align}
   a(u_h,v_h) = f(v_h) \quad \forall v_h \in V_h.
   \tag{$P_h$}
   \label{eq:vp_Vh}
\end{align}
For the adaptive finite element solution, $\cT_h$ (and thus $V_h$) is generated according to a given quantity of interest, for example the error of the solution in a chosen norm.

This study follows the $M$-uniform mesh approach \cite{Huang05a} which generates an adaptive mesh as a uniform mesh in the metric specified by a symmetric and strictly positive definite tensor $M = M(\boldx)$. 
Such a mesh is called an \emph{$M$-uniform mesh}.
Once a metric tensor $M$ has been chosen, a sequence of mesh and corresponding finite element approximation are generated in an iterative fashion.

Let $\{\cT_h^{(i)}\}$ ($i=0, 1, ...$) be an affine family of simplicial meshes on $\Omega$ with the corresponding space  $V_h^{(i)}$ of continuous, piecewise linear functions over $\cT_h^{(i)}$.
An adaptive algorithm starts with an initial mesh $\cT_h^{(0)}$.
On every mesh $\cT_h^{(i)}$  the variational problem $(P_h)$ with $V_h^{(i)}$ is solved and the obtained approximation $u_h^{(i)}$ is used to compute a new adaptive mesh for the next iteration step.
The new mesh $\cT_h^{(i+1)}$ is generated as a $M$-uniform mesh with a metric tensor $M_h^{(i)}$ defined in terms of $u_h^{(i)}$.
This yields the sequence
\[
   (\cT_h^{(0)}, V_h^{(0)}) \rightarrow u_h^{(0)} \rightarrow M_h^{(0)} \rightarrow 
   (\cT_h^{(1)}, V_h^{(1)}) \rightarrow \dots % u_h^{(1)} \rightarrow M_h^{(1)} \rightarrow \hdots 
\]
%The process is repeated until a good adaptation is achieved.
%The mesh adaptation quality is characterized by the alignment and equidistribution quality measures introduced in \cite{Huang05}; more details can be found in \cite[Sect.~4.1]{HuKaLa10}.
The conformity of the mesh to the input metric is characterized by the mesh quality measure $Q_{mesh}$ introduced in \cite{HuKaLa10}.
$Q_{mesh} \geq 1$ for any given mesh but $Q_{mesh} = 1$ if and only if the underlying mesh is $M$-uniform (see \cite[Sect.~4.1]{HuKaLa10} for more details).
In the presented numerical examples, the mesh adaptation process is repeated until the mesh is $M$-uniform within a given tolerance: $ Q_{mesh} \leq 1 + \varepsilon$ with a chosen tolerance $\varepsilon$ ($\varepsilon = 0.1$ is used in the numerical experiments throughout the paper).

Mesh generation software BAMG is used to generate new adaptive meshes (\emph{bidimensional anisotropic mesh generator} developed by F. Hecht \cite{bamg}).

\subsection{Adaptation based on \emph{a~posteriori} error estimates}
%--------------------------------------------------------------------------------
Typically, the metric tensor $M_h$ depends on the Hessian of the exact solution of the underlying problem \cite{ForPer01,Huang05a}. 
As mentioned in the introduction, it is not possible to obtain an accurate Hessian recovery from a linear finite element solution in general \cite{Kamens09PhD}, so there is no way to prove a convergence of an adaptive algorithm based on the Hessian recovery in a direct way, even if its application is successful in practical computations \cite{Dolejs98,HuaLi10,LiHua10}.

An alternative approach developed in \cite{HuKaLa10} employs an \emph{a~posteriori} error estimator for defining and computing $M_h$. 
The brief idea is as follows.

Assume that an error estimate $z_h$ is reliable in the sense that
\begin{align}
   \norm{u - u_h} \leq C \norm{z_h}.
   \label{eq:reliable}
 \end{align}
for a given norm $\norm{\cdot}$ and that it has the property
\begin{align}
   \Pi_h z_h \equiv 0 
   \label{eq:ihrh}
\end{align}
for some interpolation operator $\Pi_h$.
Then the finite element approximation error is bounded by the (explicitly computable) interpolation error of the error estimate $z_h$, viz.,
\begin{align}
  \|u - u_h\| \le C \norm{z_h} = C \norm{z_h - \Pi_h z_h}.
   \label{eq:ehrhuh}
\end{align}
Now, it is known from the interpolation theory \cite{HS01} that the interpolation error for a given function $v$ can be bounded by a term depending on the triangulation $\cT_h$ and derivatives of $v$, i.e.,
\begin{align*}
   \| v - \Pi_hv \| \leq C \,  \mathcal{E}(\cT_h,v),
\end{align*}
where $C$ is a constant independent of $\cT_h$ and $v$.
Therefore, we can rewrite \eqref{eq:ehrhuh} as
\begin{align*}
   \|u - u_h\| \le C \, \mathcal{E}(\cT_h,z_h) .
\end{align*}
In other words, up to a constant, the solution error is bounded by the interpolation error of the error estimate.
Thus, the metric tensor $M_h$ can be constructed to minimize the interpolation error of the $z_h$ and does not depend on the Hessian of the exact solution.

%-------------------------------------------------------------------------------
\subsection{Hierarchical basis \emph{a~posteriori} error estimate}
\label{sec:hbee}
One possibility to achieve the property \eqref{eq:ihrh} is to use the hierarchical basis error estimator.
The general framework can be found among others in the work of Bank and Smith \cite{BanSmi93} or Deuflhard et al. \cite{DeLeYs89}.
The approach is briefly explained as follows.

Let $e_h = u - u_h$ be the error of the linear finite element solution $u_h \in V_h$.
Then for all $v\in V$ we have
\begin{align}
   a(e_h,v) = f(v) - a(u_h,v).
   \tag{$E$}
 \end{align}

Let $\bar V_h = V_h \oplus W_h$ be a space of piecewise quadratic functions, where $W_h$ is the linear span of the quadratic edge bubble functions (a quadratic edge bubble function is defined as a product of the two linear nodal basis functions corresponding to the edge endpoints).
The error estimate $z_h$ is defined as the solution of the approximate error problem: find $z_h\in W_h$ such that
\begin{align}
   a(z_h, w_h) = f(w_h) - a(u_h,w_h) \quad \forall w_h \in W_h.
   \tag{$E_h$} \label{eq:eh}
\end{align}
The estimate $z_h$ can be viewed as a projection of the true error onto the subspace $W_h$ and relies on two assumptions:
the strengthened CBS-inequality
\[
   | a(v_h,w_h) | \leq \gamma ||| v_h ||| \, ||| w_h |||  
         \text{ with } \gamma < 1
\]
for all $v_h \in V_h$ and $w_h \in W_h$ and the saturation assumption 
\[ ||| u - u_q ||| \leq \beta ||| u - u_h ||| \quad \text{with} \quad \beta < 1, \]
i.e. the assumption that the quadratic finite element approximation $u_q$ is more accurate than the linear approximation $u_h$.

Now, if $\Pi_h$ is defined as the vertex-based, piecewise linear Lagrange interpolation then $z_h$ satisfies the condition \eqref{eq:ihrh} since the edge bubble functions vanish at vertices.
Then, if assumption \eqref{eq:reliable} holds, the finite element approximation error can be controlled by minimizing the interpolation error of $z_h$, i.e., the right-hand side in \eqref{eq:ehrhuh}.

Note that this definition of the error estimate is global and the cost of its exact solution is comparable to the cost of the quadratic finite element approximation.
To avoid the expensive exact solution in numerical computation, only a few sweeps of the symmetric Gauss-Seidel iteration are employed for the resulting linear system (until the relative difference of the old and the new approximations is under a given relative tolerance).
In numerical experiments, three to ten sweeps were enough to achieve the relative tolerance of $1\%$ and it proved to be fully sufficient for the purpose of mesh adaptation.
The computational cost of such approximation is comparable to the cost of the Hessian recovery: in the tests, the computation of HBEE was in average about two times slower than the Hessian recovery method.
Since both methods require only a fraction of the overall computational cost, using the fast HBEE approximation for mesh adaptation instead of Hessian recovery will not significantly increase the overall computational time.

\section{Poisson problem in a domain with a corner singularity}
\label{sec:cornerSingularity}
%********************************************************************
Consider the Dirichlet problem for the Poisson equation
\begin{align} 
      \begin{cases} 
         \Delta u = f      &\text{in } \Omega, \\
         u = g &\text{on } \partial\Omega,
      \end{cases} 
      \label{ex:cornerSingularity}
      \end{align}
      with $\Omega = \{x\in \R^2 : r<1, ~ 0 < \theta < \pi / \lambda \}$ for some $\frac{1}{2} < \lambda < 1$, $f=0$ and $g = r^{\lambda} \sin\left(\lambda\theta\right)$ in usual polar coordinates $(r,\theta)$.
The solution of \eqref{ex:cornerSingularity} or, more precise, the solution of the corresponding variational problem \eqref{eq:vp_H}is given by
\[
   u(r,\theta) = r^{\lambda} \sin\left(\lambda \theta\right).
\]
The gradient of the solution $u$ has a singular behaviour near the corner $(0,0)$ of the domain $\Omega$ and $u \in H^{1+\lambda-\epsilon}(\Omega)$ with arbitrarily small $\epsilon$ \cite{Johnson87}.

\subsection{Optimal mesh grading}
%----------------------------------------------------------------------
For a uniform or quasi-uniform mesh, the poor regularity of $u$ leads to the suboptimal rate of convergence of the linear finite element method:
\begin{align}
   \norm{u - u_h}_{H^1(\Omega)} &\leq C h^{\lambda - \epsilon}, \label{eq:cornerH1} \\
   \norm{u - u_h}_{L^2(\Omega)} &\leq C h^{2\lambda - \epsilon}. \label{eq:cornerL2}
\end{align}
A number of techniques were developed to regain the optimal convergence order by using specially adapted finite element spaces, see \cite[Sect.~4.2]{Apel99} and the references therein for more details.
In particular, it is possible to regain the optimal convergence rate by using the standard linear finite element space with proper mesh grading around the corner of the domain.

The basic idea is described among others in \cite{Johnson87} and is as follows.
The linear finite element error is bounded by the linear interpolation error
\[
   \abs{u - u_h}^2_{H^1(\Omega)} 
      \leq \abs{u - \Pi_h u}^2_{H^1(\Omega)}
      \leq C \sum_{K\in \cT_h} h^2_K \abs{u}^2_{H^2(K)}.
\]
Thus, in order to obtain the optimal convergence rate, the mesh size $h_k$ has to be balanced with the size of $\abs{u}_{H^2(K)}$ .
Now, if mesh elements in a neighbourhood around the corner $(0,0)$ have the size
\begin{align*}
   h_K = C h r_K^{1-\mu},
   %\label{eq:cornerDistribution}
\end{align*}
where $h$ is the mesh size far away from the corner and $r_K$ is the distance of an element $K$ to the corner, then the linear finite element error is
\begin{align*}
   \norm{u - u_h}_{H^1(\Omega)} &\leq C h,\\
   \norm{u - u_h}_{L^2(\Omega)} &\leq C h^2,
\end{align*}
provided that $\mu < \lambda$ \cite{Apel99,Johnson87}. Note, that such meshes still have $N \sim h^{-2}$ number of elements, thus significantly increasing the accuracy of the approximation for a given number of mesh elements.

In this example, we follow the approach described in the Sect.~\ref{sec:metricBasedAdaptation} and employ the metric tensor developed in \cite{HuKaLa10} for a boundary  value problem of a second-order elliptic PDE.
In two dimensions and for the $L^2$-norm of the error, the metric tensor $M_{HB}$ based on the HB error estimator is given element-wise as
\begin{align*}
  M_{HB,K} = 
   \det\left(I + \frac{1}{\alpha_h} \abs{H_K(z_{h})}\right)^{-\frac{1}{6}}
   \left[I + \frac{1}{\alpha_h} \abs{H_K(z_{h})}\right],
   %\label{eq:madapt}
\end{align*}
where $H_K(z_h)$ is the Hessian of the quadratic hierarchical basis error estimate $z_h$ on element $K$ and $\alpha_h$ is a regularization parameter to ensure that the metric tensor is strictly positive definite.
Also, $\alpha_h$ can be used as a adaptation intensity control: if $\alpha_h \rightarrow \infty$, the mesh becomes uniform; if $\alpha_h \rightarrow 0$, the mesh becomes more adaptive.
Usually, $\alpha_h$ is chosen so that about half of the mesh elements are concentrated in regions where $\det(M)$ is large (see \cite{HuKaLa10} for more details on the choice of $\alpha_h$).

Note, that $M_{HB}$ is completely \emph{a~posteriori} and does not require any \emph{a~priori} information on the solution.

\subsection{Numerical example}
%--------------------------------------------------------------------------------
Consider the Dirichlet boundary problem \eqref{ex:cornerSingularity} with $\lambda = 4/7$.
In two dimensions, the number of elements $N$ of a quasi-uniform mesh is proportional to $h^{-2}$, so there is a factor of $-1/2$ in the convergence rate if it is given in terms of the number of mesh elements. 
Thus, according to \eqref{eq:cornerH1} and \eqref{eq:cornerL2}, the expected orders of convergence in terms of the number of mesh elements for quasi-uniform meshes should be at most $-2/7$ (i.e. $-\lambda/2$) for the $H^1$ and $-4/7$ (i.e. $-\lambda$) for the $L^2$ norms of the error.
For the optimally graded meshes we should expect orders $-0.5$ and $-1$, respectively.

\begin{figure}[t] \centering
  \subfloat[Quasi-uniform, 1214 triangles.]{  % 7 525
      \includegraphics[width=0.23\textwidth,clip]{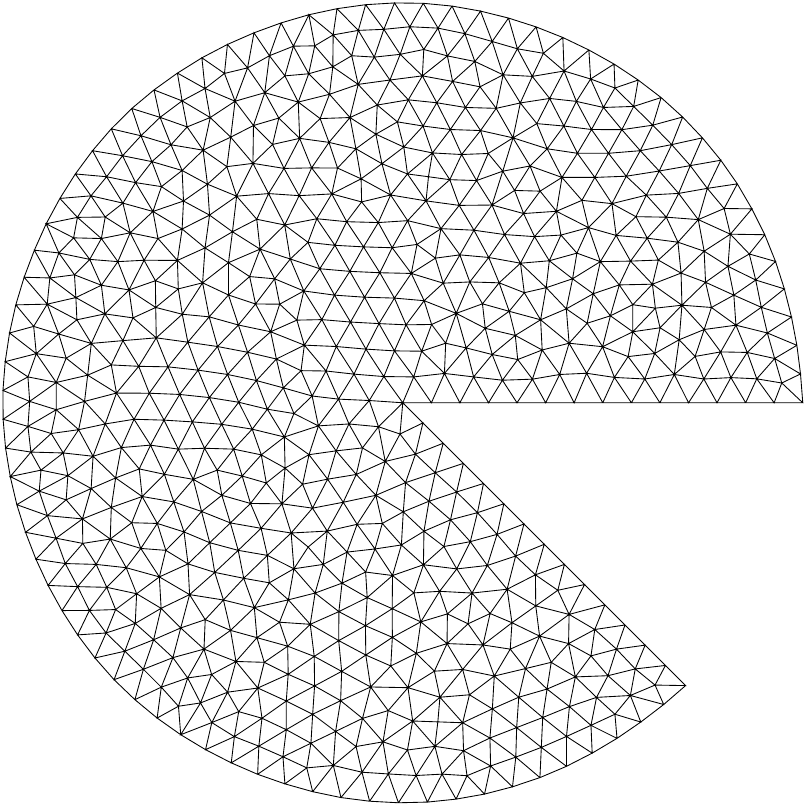}
      \label{fig:cornerUniform}
   } 
   \subfloat[$M_{HB}$, 1225 triangles.]{   % 5 500
      \includegraphics[width=0.23\textwidth,clip]{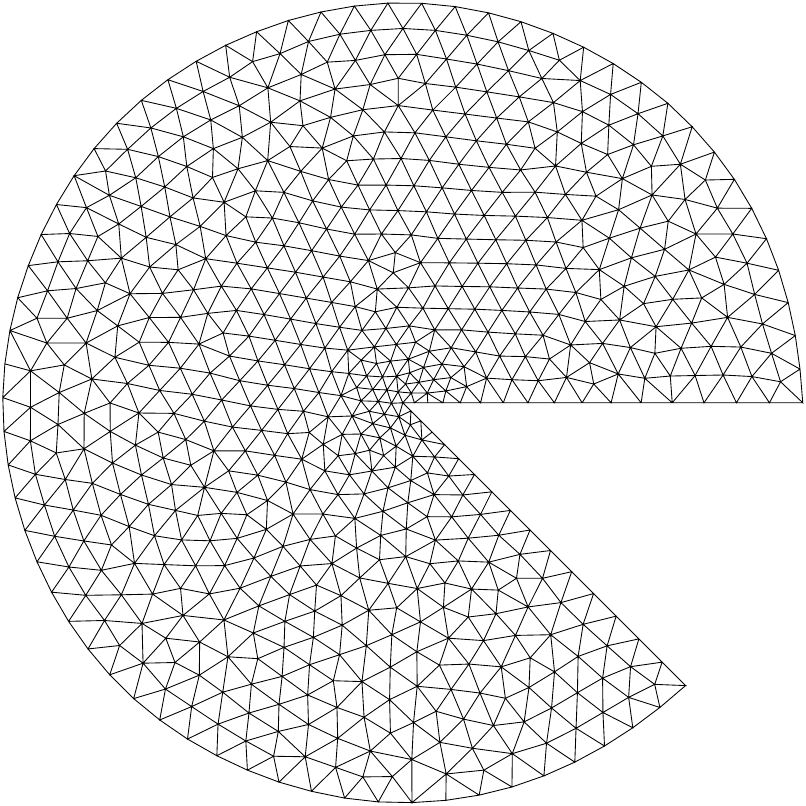}
      \label{fig:cornerAdapt}
   }
   \caption{Corner singularity and optimal mesh grading: mesh examples.}
   \label{fig:cornerMeshes}
\end{figure}
\begin{figure}[t] \centering
   \subfloat[Quasi-uniform: 165\,936 triangles.]{ % 7 88000
      \includegraphics[width=0.23\textwidth,clip]{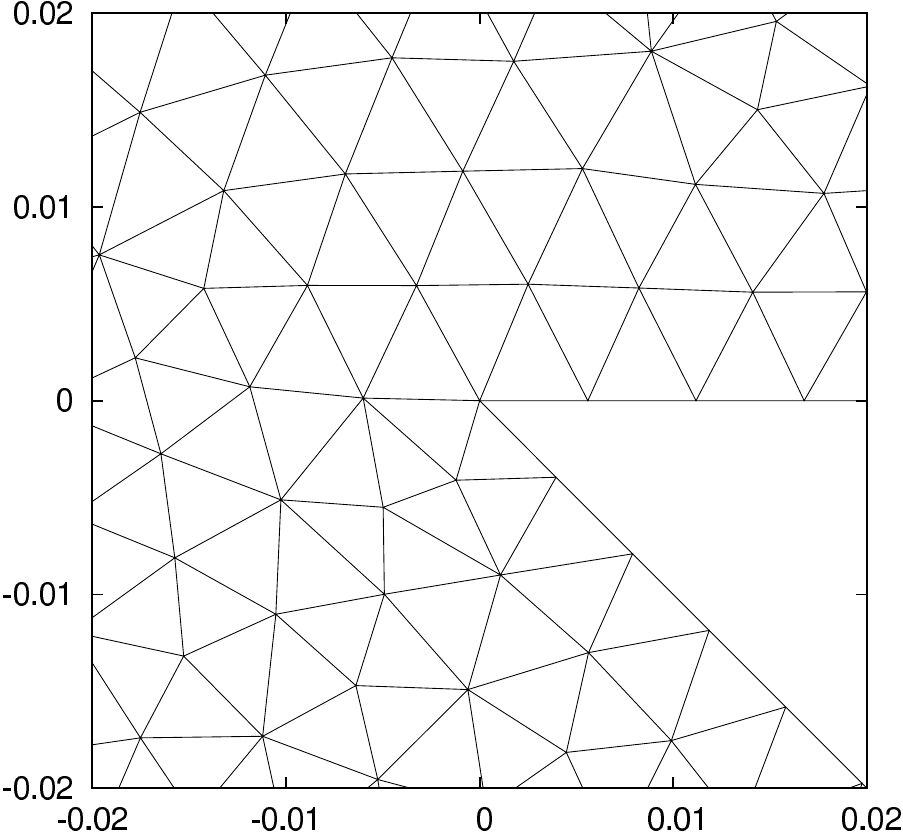}
      \label{fig:cornerUniformZoom}
   } 
   \subfloat[$M_{HB}$, 166\,125 triangles.]{ % 5 80000
      \includegraphics[width=0.23\textwidth,clip]{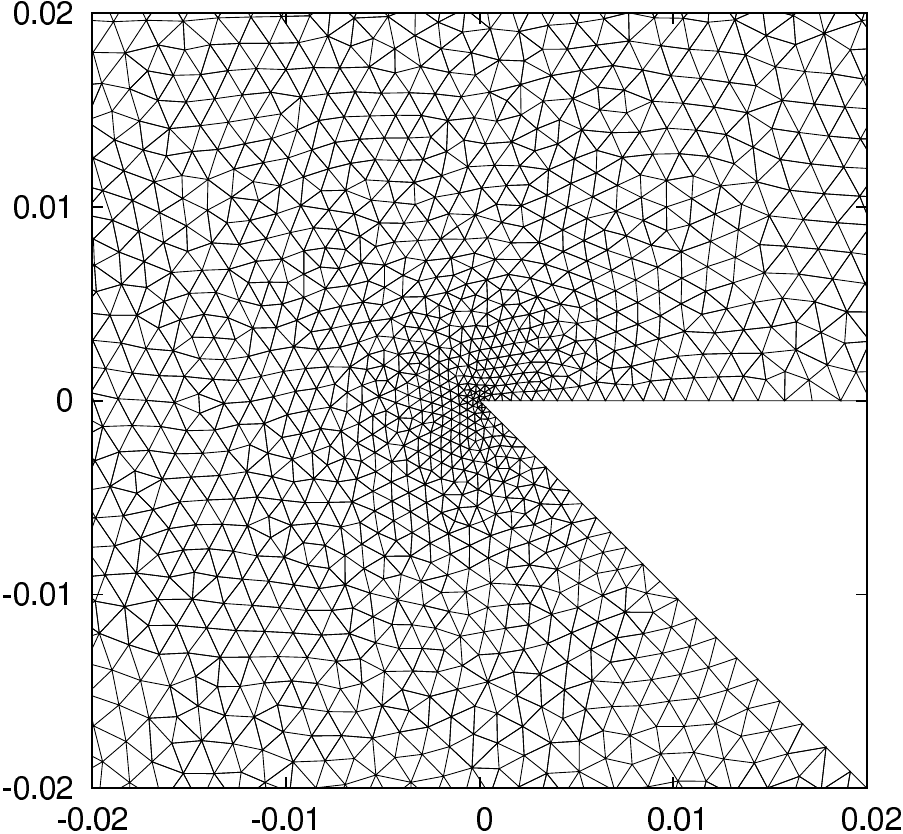}
      \label{fig:cornerAdaptZoom}
   }
   \caption{Corner singularity and optimal mesh grading: close-ups at $(0,0)$.}
   \label{fig:cornerMeshesAdapt}
\end{figure}
\begin{figure}[t] \centering
   \includegraphics[width=0.48\textwidth,clip]{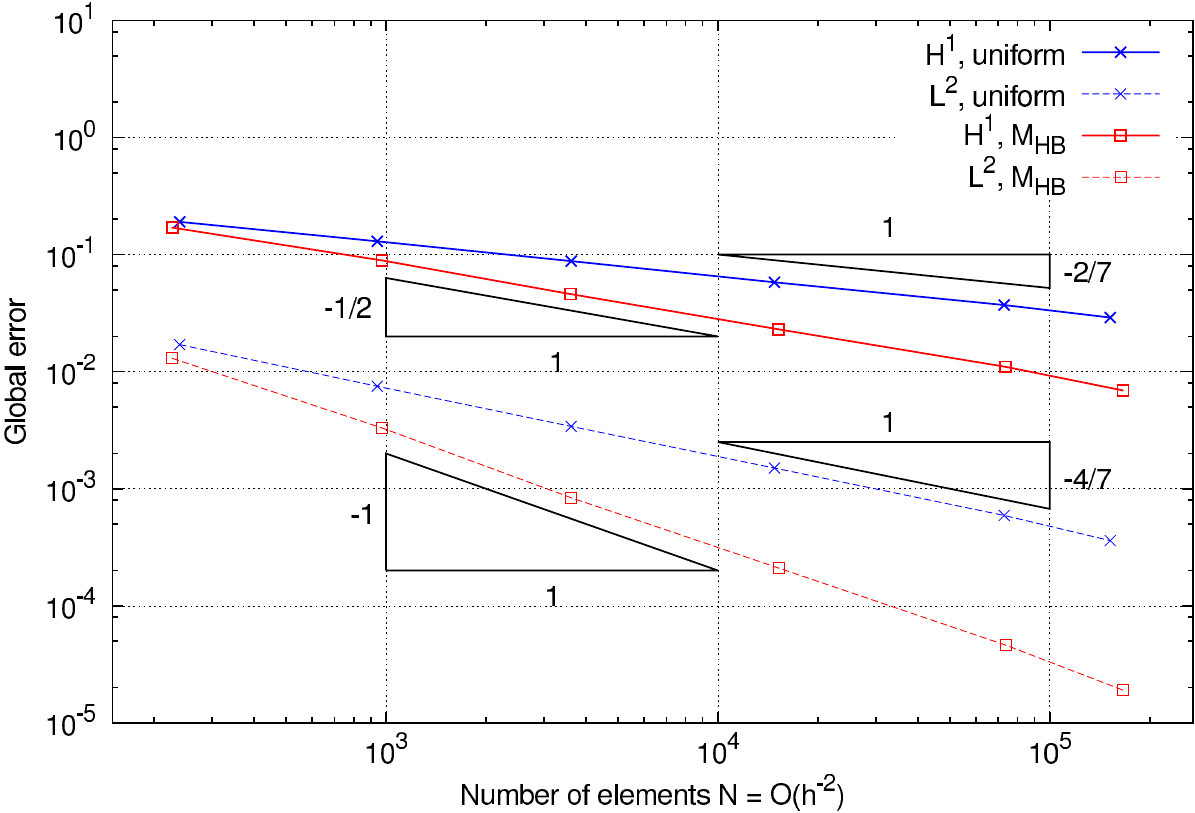}
   \caption{Corner singularity and optimal mesh grading: global error 
      for quasi-uniform and adaptive meshes against the number of elements.}
   \label{fig:cornerConv}
\end{figure}

Examples of uniform and adaptive meshes as well as close-up views near the corner are given in Figs.~\ref{fig:cornerMeshes} and \ref{fig:cornerMeshesAdapt}.
For the adaptive mesh, the concentration of mesh points near the singularity can be clearly observed (Figs.~\ref{fig:cornerAdapt} and \ref{fig:cornerAdaptZoom}).
We can also observe the gradual change of mesh elements size in the close-up of the adaptive mesh (Fig.~\ref{fig:cornerAdaptZoom}); in particular, along the boundary edges from the domain corner to the outer boundary.
The mesh grading is rather moderate, but it turns out to be enough to achieve the optimal convergence order, as can be observed in the corresponding convergence plot (Fig.~\ref{fig:cornerConv}).

Figure~\ref{fig:cornerConv} shows the $H^1$ and $L^2$ norms of the global linear finite element solution error on uniform and adaptive meshes against the number of mesh elements.
As expected, the convergence order of the finite element solution on uniform meshes is only $-2/7$ and $-4/7$ for the $H^1$ and the $L^2$ norms of the error, respectively.
Contrary, the convergence order of the error of the finite element solution obtained by means of the metric tensor $M_{HB}$ is $-1/2$ and $-1$. 

Thus, the considered method is able to achieve the optimal mesh grading in this example relying on only on the \emph{a~posteriori} information provided by the hierarchical basis error estimator.

\section{Variational problems}
\label{sec:variationalProblem}
%--------------------------------------------------------------------
While it has attracted considerable attention from many researchers and been successfully applied to the numerical solution of PDEs, anisotropic mesh adaptation has rarely been employed for variational problems, especially when combined with \emph{a~posteriori} error estimates.
Recently, Huang and Li \cite{HuaLi10} developed a metric tensor for the adaptive finite element solution of variational problems.
In the anisotropic case, it is semi-\emph{a~posteriori}: it involves residual and edge jumps, both dependent on the computed solution, and the Hessian of the exact solution.
In \cite{HuKaLi10}, this result was improved to provide a metric tensor for variational problems based on the HBEE and the underlying variational formulation.
The new metric tensor is \emph{a~posteriori} in the sense that it is based solely on residual, edge jumps, and HBEE.
This is in contrast to most previous work where $M$ depends on the Hessian of the exact solution and is semi-\emph{a~posteriori} or completely \emph{a~priori}; e.g., see \cite{BoGHLS97,BoGeMo97,CaHeMP97,Frey08,Huang05a,HuaLi10}.

\subsection{General variational problem and the anisotropic metric tensor}
%--- problem definition
Consider a general functional of the form 
\begin{align*}
   I[v] = \int_{\Omega} F(\boldx,v,\nabla v) d\boldx ,
      \quad \forall v \in V_{g}
\end{align*}
where $F(\cdot, \cdot, \cdot)$ is a given smooth function, $\Omega \subset \R^{d}$ ($d = 1,2,3$) is the physical domain and $V_{g}$ is a properly selected set of functions satisfying the Dirichlet boundary condition 
\begin{align*}
   v(\boldx) = g(\boldx) \quad \forall \boldx \in \partial \Omega
\end{align*}
for a given function $g$.

The corresponding variational problem is to find a minimizer $u \in V_{g}$ such that
\begin{align*}
   I[u]=\min_{v\in V_{g}}I[v].
\end{align*}
A necessary condition for $u$ to be a minimizer is that the first variation of the functional vanishes.
This leads to the Galerkin formulation 
\begin{align}
   \delta I [u,v] 
   \equiv \int_{\Omega} \left( F_{u}(\boldx,u,\nabla u) \; v 
         + F_{\nabla u}(\boldx,u,\nabla u)\cdot \nabla v\right) d\boldx 
   = 0
   \label{eq:galerkin1}
\end{align}
for all $v \in V_{0}$, where $V_{0} = V_{g}$ with $g = 0$ and $F_{u}$ and $F_{\nabla u}$ are the partial derivatives of $F$ with respect to $u$ and $\nabla u$, respectively.

Given a triangulation $\cT_h$ for $\Omega$ and the associated linear finite element space $V_{g,h} \subset V_g$, the finite element solution $u_h$ can be found by solving the corresponding Galerkin formulation: find $u_h \in V_{g,h} $ such that
\begin{align*}
   %\delta I [u_h, v_h] 
      \int_\Omega \left(F_u(\boldx,u_h,\nabla u_h) \; v_h 
         + F_{\nabla u}(\boldx,u_h,\nabla u_h) \cdot \nabla v_h \right) d\boldx
      = 0
\end{align*}
for all $v_h \in V_{0,h}$.

The \emph{a~posteriori} metric tensor $M_{HB,K}$ for general variational problems developed in \cite{HuKaLi10} for the error measured in the $H^1$-semi-norm is given element-wise by
\begin{align*}
  M_{HB,K} &=  \left(1 + \frac{1}{\alpha_h \abs{K}}
  \left( \abs{K}^{1/2} \norm{r_h}_{L^2(K)} \right. \right. \notag\\
  & \qquad\quad \left.\left. + \sum_{\gamma\in\partial K} \abs{\gamma}^{1/2}
            \norm{R_h}_{L^2(\gamma)}\right)\right)^{\frac{1}{2}} 
      \notag \\
      & \quad \times 
      \det \left(I + \frac{1}{\alpha_h}\abs{H_K(z_h)}\right)^{-\frac{1}{4}}
      \left[I + \frac{1}{\alpha_h}\abs{H_K(z_h)}\right]
   % \label{metric-an1}
\end{align*}
with the residual
\begin{align*}
   r_{h}(\boldx)  
      = F_{u}(\boldx)-\nabla \cdot F_{\nabla u}(\boldx) 
         \qquad \boldx \in K \quad \forall K \in \cT_h
\end{align*}
and the edge jump
\begin{align*}
R_{h}(\boldx) = 
   (F_{\nabla u}(\boldx)\cdot \boldsymbol{n}_\gamma)|_{K}
      +(F_{\nabla u}(\boldx)\cdot\boldsymbol{n}_\gamma)|_{K'}
\end{align*}
on $\gamma \in \partial \cT_h\backslash \partial \Omega$ ($R_h = 0$ if  $\gamma \in \partial \Omega$).
% \begin{align*}
% R_{h}(\boldx) = 
%    \begin{cases}
%    (F_{\nabla u}(\boldx)\cdot \boldsymbol{n}_\gamma)|_{K}
%       +(F_{\nabla u}(\boldx)\cdot\boldsymbol{n}_\gamma)|_{K'}
%       & \text{on } \gamma \in \partial \cT_h\backslash \partial \Omega,\\ 
%       0 & \text{on } \gamma \in \partial \Omega.
%    \end{cases}
% \end{align*}
As in the previous section, $\alpha_h$ is a regularization parameter to ensure that the metric tensor is strictly positive definite (see \cite{HuKaLa10} for more details on the choice of $\alpha_h$).

Note, that $\delta I[u,v]$ in \eqref{eq:galerkin1} is linear in $v$ but is nonlinear in $u$ in general. 
Thus, a modification of the error problem \eqref{eq:eh} for $z_h$ in Sect.~\ref{sec:hbee} is required.
For this purpose, denote by $a_h(u_h;\cdot,\cdot)$ a bilinear form resulting from a linearization of $\delta I[\cdot,\cdot]$ about $u_h$ with respect to the first argument.
The error estimate $z_h$ is then defined as the solution of the approximate \emph{linear} error problem: find $z_h\in W_h$ such that
\begin{align*}
     a_h(u_h;z_h,w_h) =  - \delta I[u_h,w_h]
\end{align*}
for all $w_h \in W_h$ \cite{BanSmi93}.

\subsection{Numerical example}
%--- Adaptive linear finite element approximation ------------------------
\begin{figure}[t] \centering
   \includegraphics[width=0.30\textwidth,clip]{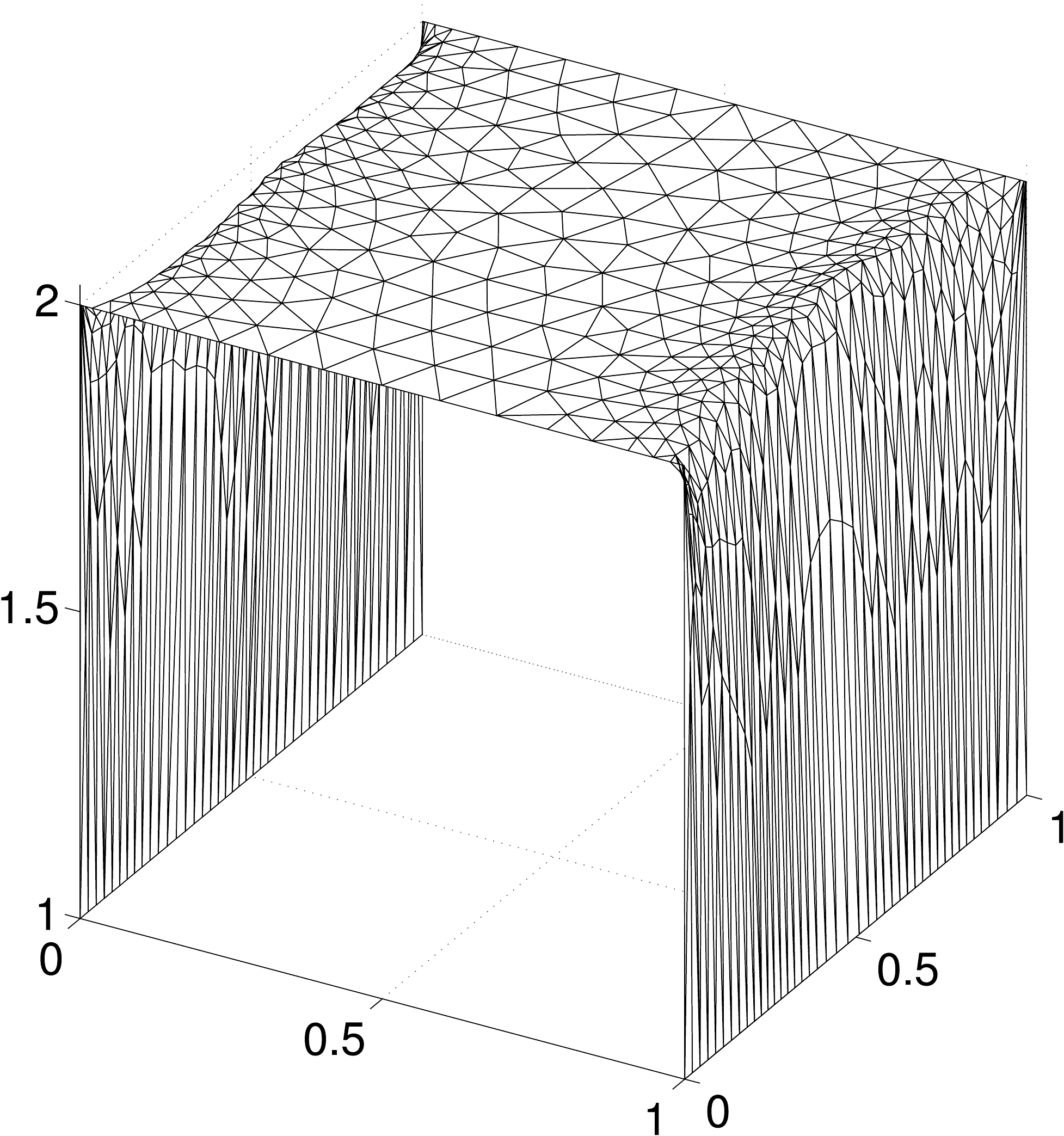}
   \caption{Variational problem: numerical solution.}
   \label{fig:plot:2}
\end{figure}

Consider an anisotropic variational problem defined by the non-quadratic functional
\begin{align*}
   I[u] = \int_\Omega \left[ \left( 1 + \abs{\nabla u}^2 \right)^{3/4}
      + 1000 u_y^2 \right] d\boldx
\end{align*}
with $\Omega = (0,1) \times (0,1)$ and the boundary condition 
\begin{align*}
   \begin{cases}
      u = 1  & \text{on } x=0 \text{ or } x=1, \\
      u = 2  & \text{on } y=0 \text{ or } y=1.
   \end{cases}
\end{align*}
This example is discussed in \cite{HuKaLi10,HuaLi10} and is originally taken from \cite{Bildhauer03}; the analytical solution is not available, but a computed solution in Fig.~\ref{fig:plot:2} shows that the mesh adaptation challenge for this example is the resolution of the sharp boundary layers near $x = 0$ and $x = 1$.

% non-linear functional example: qls<->hb: mesh examples
% beta = 0.8; 510 (qls), 500 (hb)
\begin{figure}[t] \centering
   \subfloat[Metric based on residual, edge jumps, and Hessian recovery:
         1160 triangles, maximum aspect ratio 15.]{
      \includegraphics[width=0.23\textwidth,clip]{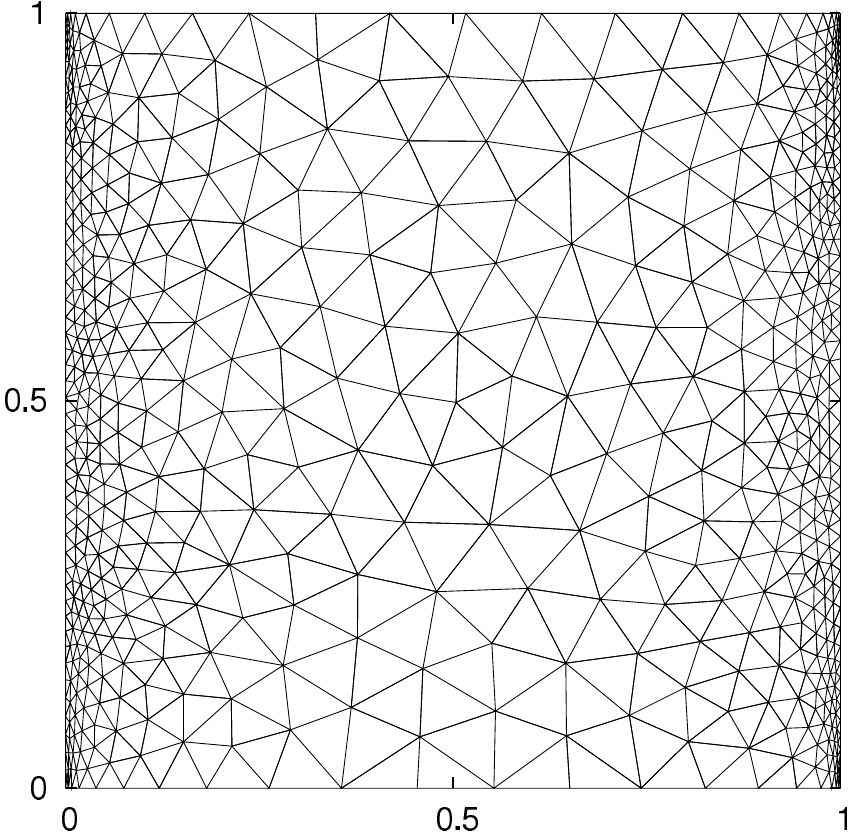}
      \includegraphics[width=0.238\textwidth,clip]{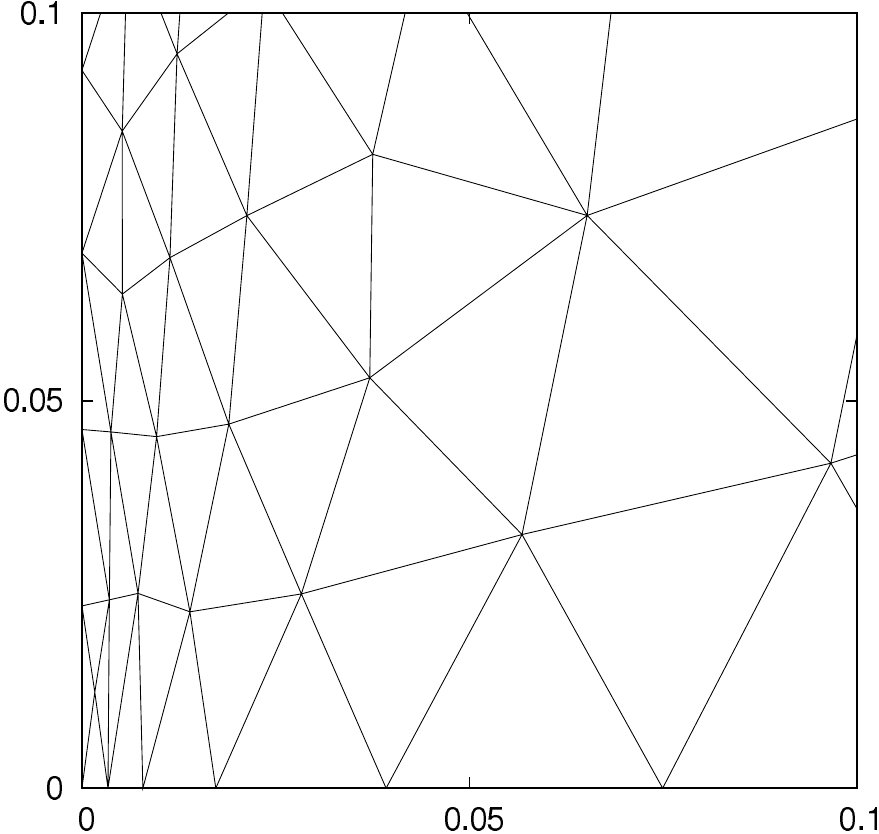}
      \label{fig:nonLinearExample:qls}
   }\\
   \subfloat[Metric based on residual, edge jumps, and HBEE:
            1143 triangles, maximum aspect ratio 51.]{
      \includegraphics[width=0.23\textwidth,clip]{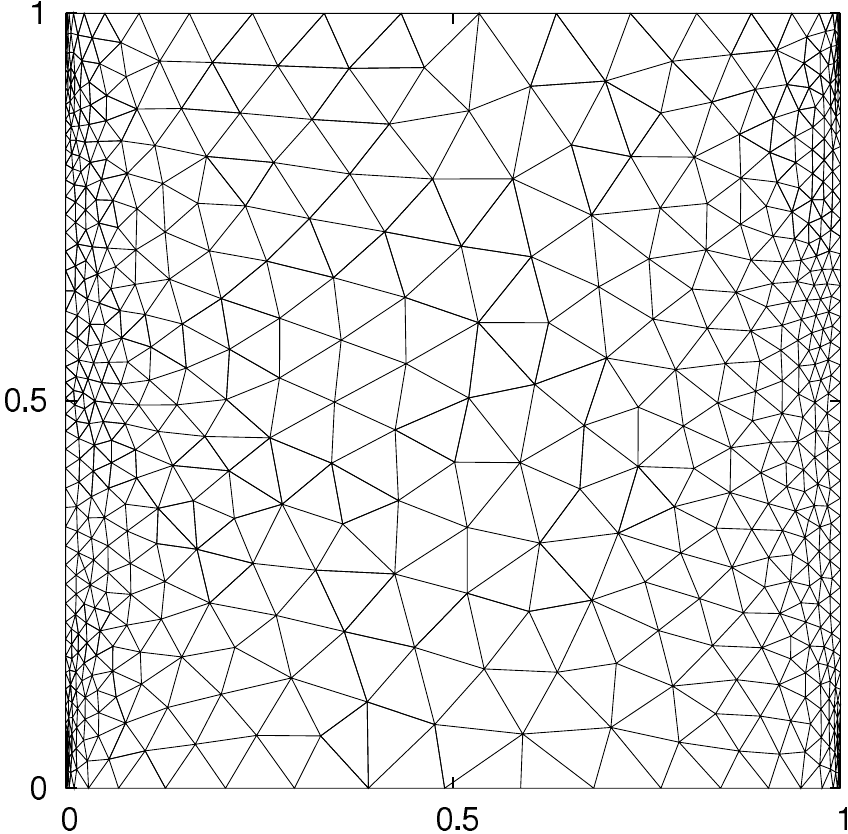}
      \includegraphics[width=0.238\textwidth,clip]{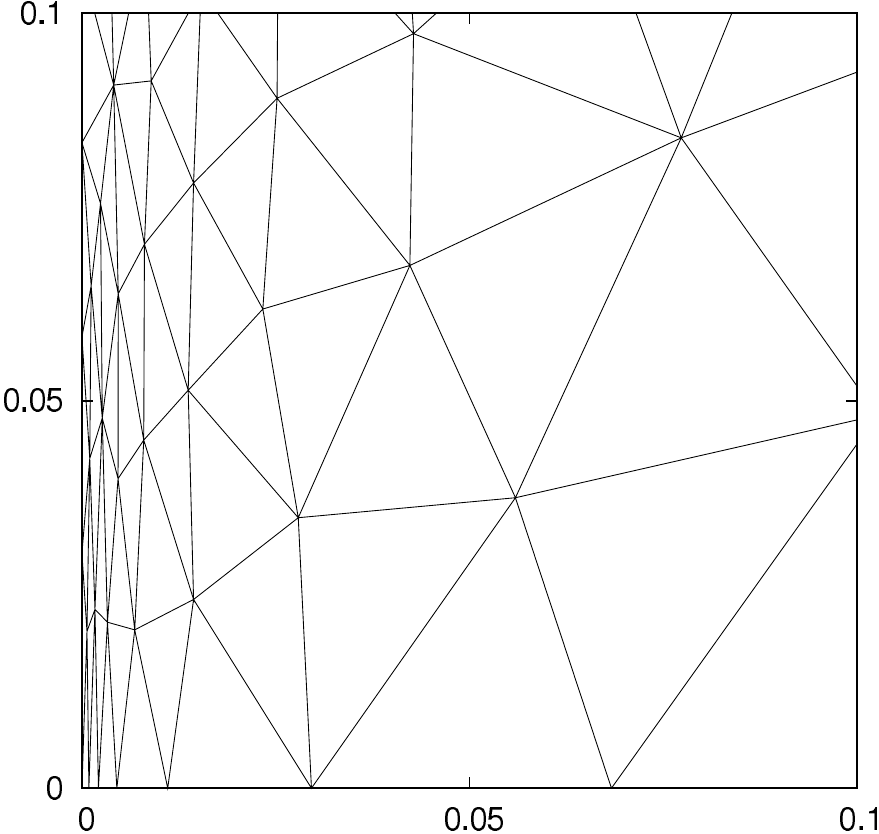}
      \label{fig:nonLinearExample:hb}
   }
     \caption{Variational problem: adaptive meshes and close-up views at (0,0).}
     \label{fig:nonLinearExample}
\end{figure}

Adaptive meshes obtained by means of Hessian recovery\footnote{Quadratic least squares fitting Hessian recovery \cite{ZhaNag05}, which seems to be the most robust and reliable Hessian recovery method \cite{Kamens09PhD,VaMDDG07}.}
and HBEE are given in Fig.~\ref{fig:nonLinearExample}.
The both methods have correct mesh concentration and provide good alignment with the boundary layers.
Anisotropic meshes are comparable, although mesh elements near the boundary layer in the HBEE-based adaptive mesh have a larger aspect ratio\footnote{In this paper, aspect ratio is defined as the longest edge divided by the shortest altitude. An equilateral triangle has an aspect ratio of $2 / \sqrt{3} \approx 1.15$.} than elements of the mesh obtained by means of the Hessian recovery.
This could be due to the smoothing nature of the Hessian recovery: usually, it operates on a larger patch, thus introducing an additional smoothing effect, which affects the grading of the elements' size and orientation.
The global hierarchical basis error estimator does not have this handicap and, in this example, the mesh obtained by means of HBEE has a better refinement in the orthogonal direction along the steep boundary layers.

\section{Anisotropic diffusion and the DMP}
\label{sec:anisotopicDiffusion}
%--------------------------------------------------------------------
Anisotropic diffusion problems arise in various areas of science and engineering, for example image processing, plasma physics, or petroleum engineering.
Standard numerical methods can produce spurious oscillations when they are used to solve these problems. 
A common approach to avoid this difficulty is to design a proper numerical scheme or a mesh so that the numerical solution satisfies the discrete counterpart (DMP) of the maximum principle satisfied by the continuous solution.
A well known condition for the DMP satisfaction by the linear finite element solution of isotropic diffusion problems is the non-obtuse angle condition that requires the dihedral angles of mesh elements to be non-obtuse \cite{CiaRav73}. 
In \cite{LiHua10}, a generalization of the condition, the so-called anisotropic non-obtuse angle condition, was introduced for the finite element solution of heterogeneous anisotropic diffusion problems. 
The new condition is essentially the same as the existing one except that the dihedral angles are measured in a metric depending on the diffusion matrix of the underlying problem, i.e. the mesh is aligned with the major diffusion directions.
Based on the new condition, a metric tensor for anisotropic mesh adaptation was developed, which combines the satisfaction of the DMP with mesh adaptivity: in two dimensions and for the error measured in the $H^1$-semi-norm, it is given element-wise by
\begin{equation}
  M_{DMP,K} = \left ( 1+ \frac{1}{\alpha_h}  B_{H,K} \right )^{\frac{1}{2}}
\mbox{det} \left ( \mathbb{D}_K \right )^{\frac{1}{2}} \mathbb{D}_K^{-1},
    \label{eq:mDMPAdapt}
\end{equation}
where
\begin{equation}
  B_{H,K} = \det \left ( \mathbb{D}_K \right )^{-\frac{1}{2}}
   \norm{\mathbb{D}_K^{-1} }\cdot  \frac{1}{\abs{K}} \int_K \norm{\mathbb{D}_K \abs{H(u)}}^2 d \boldx
 \label{eq:B}
\end{equation}
and $\alpha_h$ is a regularization parameter to ensure that the metric tensor is strictly positive definite (see \cite{HuKaLa10} for more details on the choice of $\alpha_h$).
The diffusion matrix $\mathbb{D}$ in \eqref{eq:mDMPAdapt} provides the correct mesh alignment and thus the satisfaction of the DMP whereas the Hessian of the exact solution in \eqref{eq:B} provides adaptivity with the solution profile.

As noted in Sect.~\ref{sec:metricBasedAdaptation}, if the finite element solution error can be bounded by the interpolation error of the hierarchical basis error estimator, the Hessian of the exact solution in \eqref{eq:B} can be replaced by the Hessian of the HBEE $z_h$:
\begin{equation}
  B_{HB,K} = \det \left ( \mathbb{D}_K \right )^{-\frac{1}{2}}
  \norm{\mathbb{D}_K^{-1} }\cdot  \frac{1}{\abs{K}} \int_K  
      \norm{\mathbb{D}_K \abs{H(z_h)}}^2 d \boldx.
 \label{eq:Bhbee}
\end{equation}
With this choice, the mesh will still satisfy the DMP because the mesh is still aligned with the major diffusion directions, but this time the mesh density is determined by the \emph{a~posteriori} error estimate $z_h$.

\subsection{Numerical example}
%--- Adaptive linear finite element approximation ------------------------
% figure ex1-domain
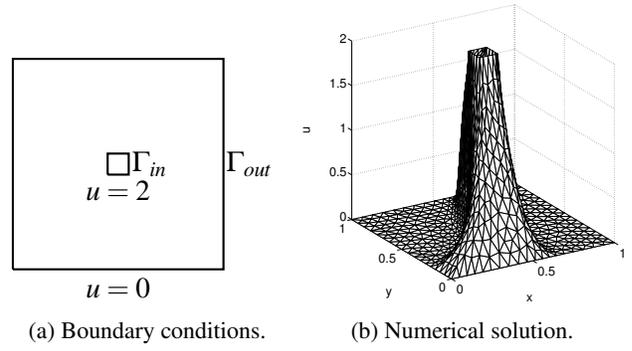
\begin{figure}[t] \centering
   \subfloat[Boundary conditions.]{
      \begin{tikzpicture}[scale = 0.7]
         \draw [thick] (0,0) -- (0,4) -- (4, 4) -- (4, 0) -- (0, 0);
         \draw [thick] (1.8,1.8) -- (1.8,2.2) -- (2.2, 2.2) -- (2.2, 1.8) -- (1.8, 1.8);
         \draw (2, -0.35) node {$ u = 0 $};
         \draw (4.5, 2) node {$ \Gamma_{out} $};
         \draw (2, 1.5) node {$ u = 2 $};
         \draw (2.6, 2) node {$ \Gamma_{in} $};
      \end{tikzpicture}
      \label{fig:aniso1_domain}
   } 
   \subfloat[Numerical solution.]{
      \includegraphics[width=0.25\textwidth,clip]{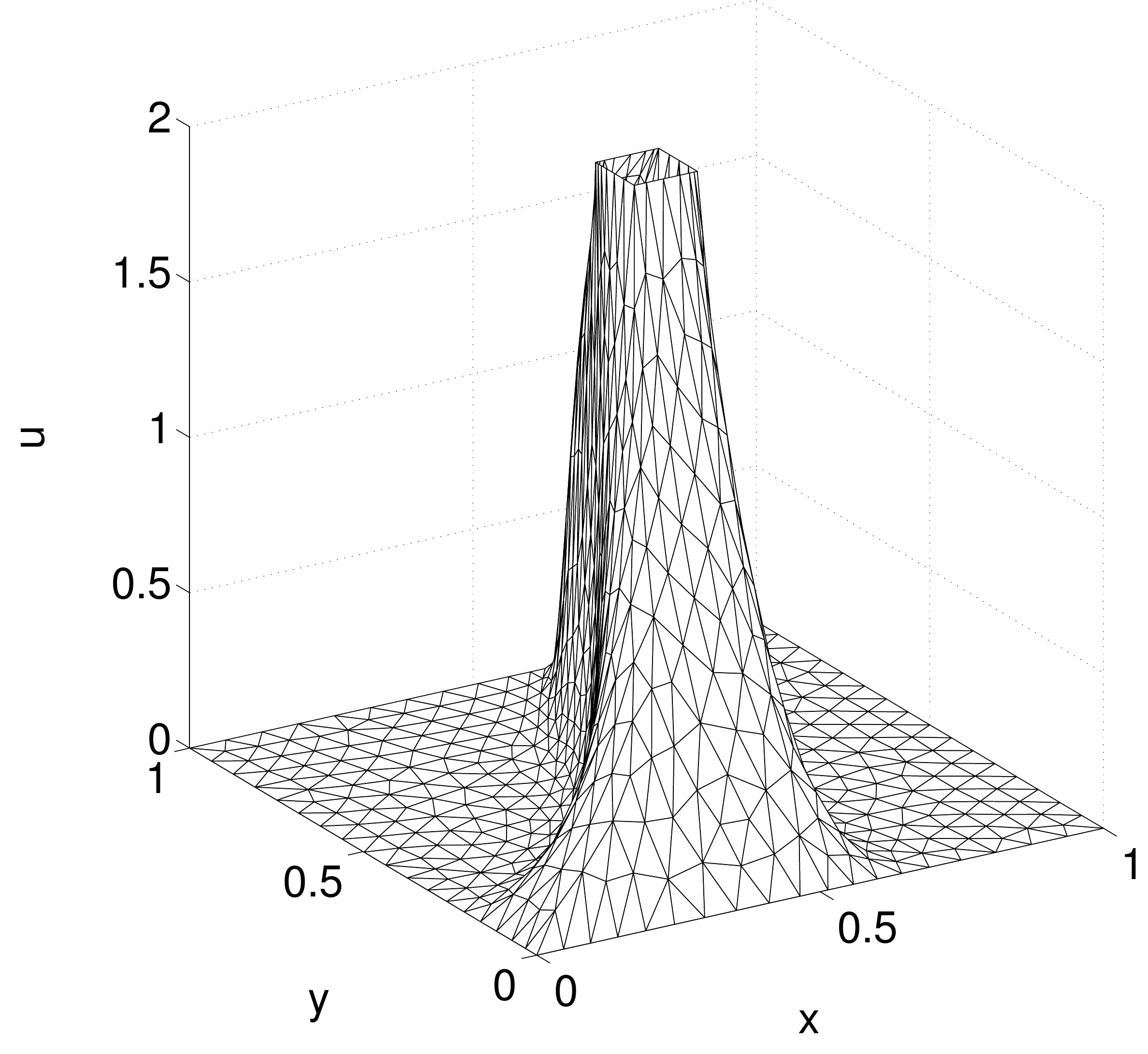}
      \label{fig:aniso1_nsol}
   }
   \caption{Anisotropic diffusion example.}
   \label{ex1-domain}
\end{figure}

Consider the BVP discussed in \cite{LiHua10}: 
\begin{align*} 
   \begin{cases} 
     \nabla \cdot (\mathbb{D} \, \nabla u) = f      &\text{in } \Omega, \\
      u = g &\text{on } \partial\Omega,
   \end{cases} 
   %\label{ex:anisodiffusion}
\end{align*}
with
\[
f\equiv 0,\quad \Omega = [0,1]^2\backslash\left [\frac{4}{9},\frac{5}{9}\right ]^2,
\quad g = 0 \mbox{ on } \Gamma_{out}, \;  g = 2 \mbox{ on } \Gamma_{in},
\]
where $\Gamma_{out}$ and $\Gamma_{in}$ are the outer and inner boundaries of $\Omega$, respectively (Fig.~\ref{fig:aniso1_domain}).
The diffusion matrix is given by 
\[
\mathbb{D} = 
   \begin{bmatrix}
      \cos\theta   & - \sin\theta \\ 
      \sin \theta  &  \cos\theta 
   \end{bmatrix}
   \begin{bmatrix}
     1000   & 0 \\ 
     0      & 1
   \end{bmatrix}
   \begin{bmatrix}
      \cos\theta     & \sin\theta \\
      - \sin\theta   & \cos\theta
   \end{bmatrix}
   % \quad \theta= \pi / 4, 
\]
with the angle of the primary diffusion direction
\[
   \theta = \pi \sin x \cos y
\] 
(the primary diffusion direction is parallel to the first eigenvector of $\mathbb{D}$).
 
This example satisfies the maximum principle and the solution stays between 0 and 2 and has sharp jumps near the inner boundary. The analytical solution is not available, but the numerical solution is provided in Fig.~\ref{fig:aniso1_nsol}.  
The goal is to produce a numerical solution which also satisfies DMP, i.e. stays between 0 and 2, and has a good adaptation.

\begin{figure}[t] \centering
  \subfloat[Isotropic: 4170 triangles, $\min u_h \approx -5.9 \times 10^{-2}$,
   maximum aspect ratio 2.7.]{
      \includegraphics[width=0.23\textwidth,clip]{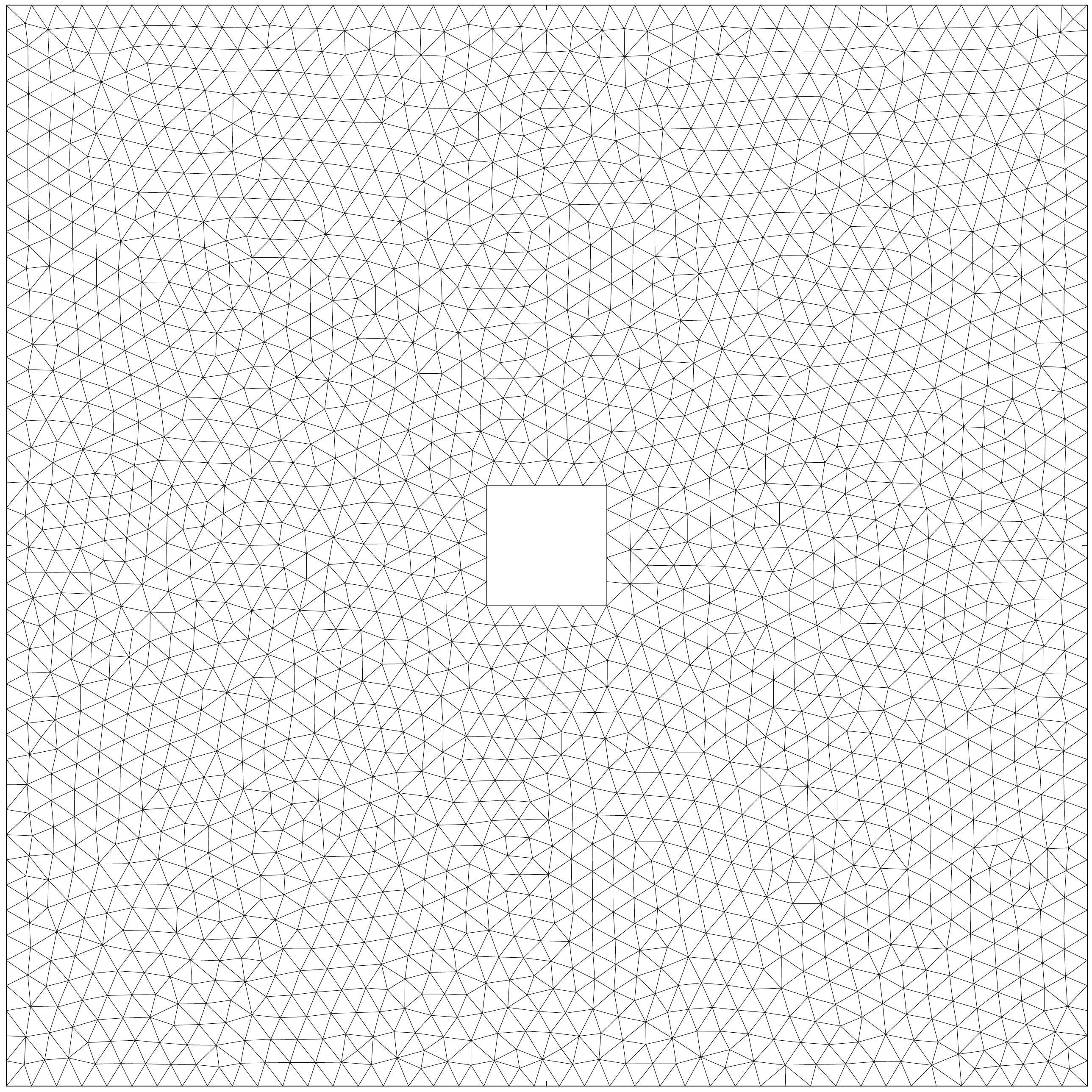}
      \includegraphics[width=0.23\textwidth,clip]{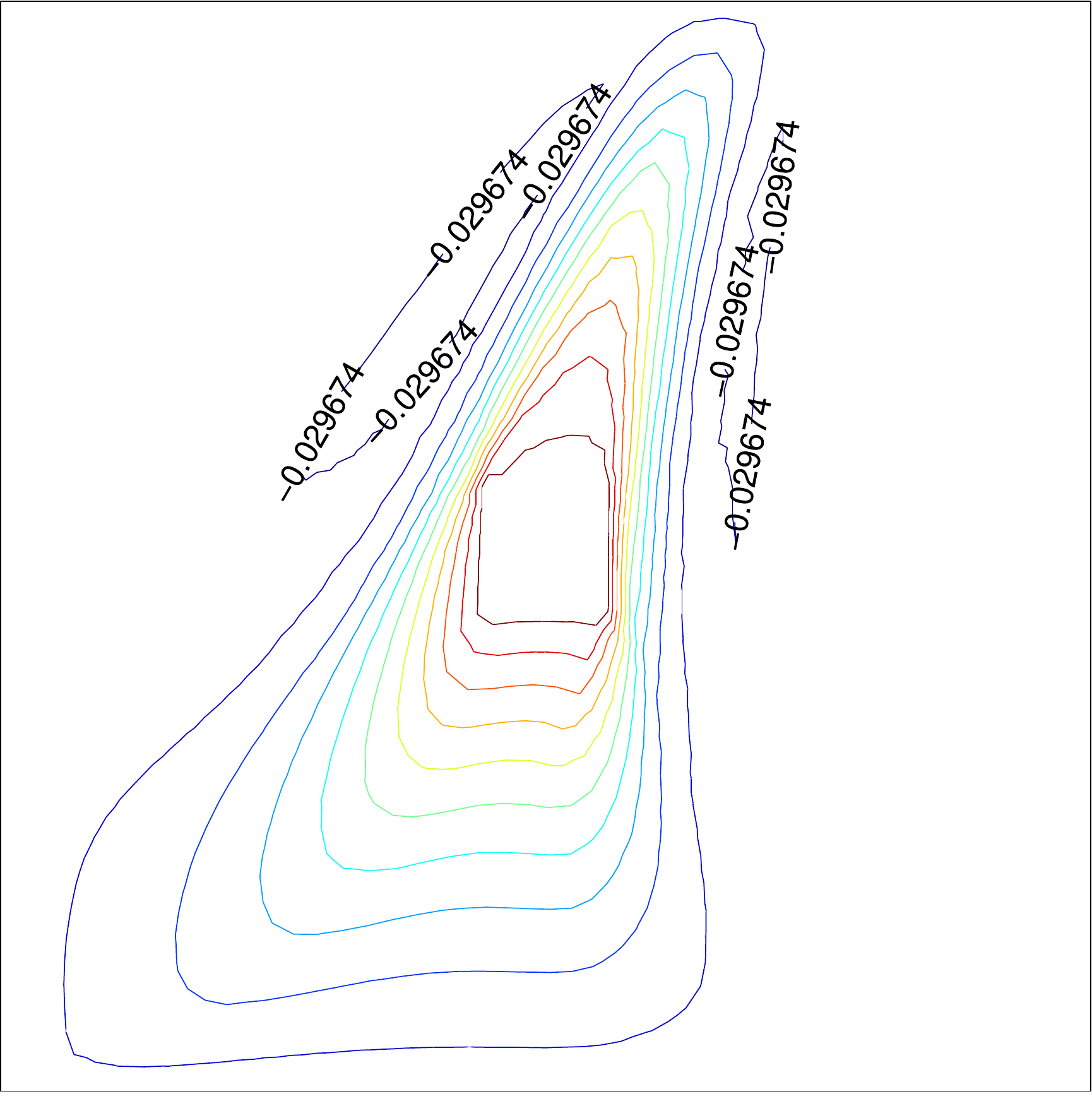}
      \label{fig:aniso1hb1}
   } \\
   \subfloat[$M_{HB}$: 4353 triangles, $\min u_h \approx -3.2 \times 10^{-4}$,
   maximum aspect ratio 37.2.]{
      \includegraphics[width=0.23\textwidth,clip]{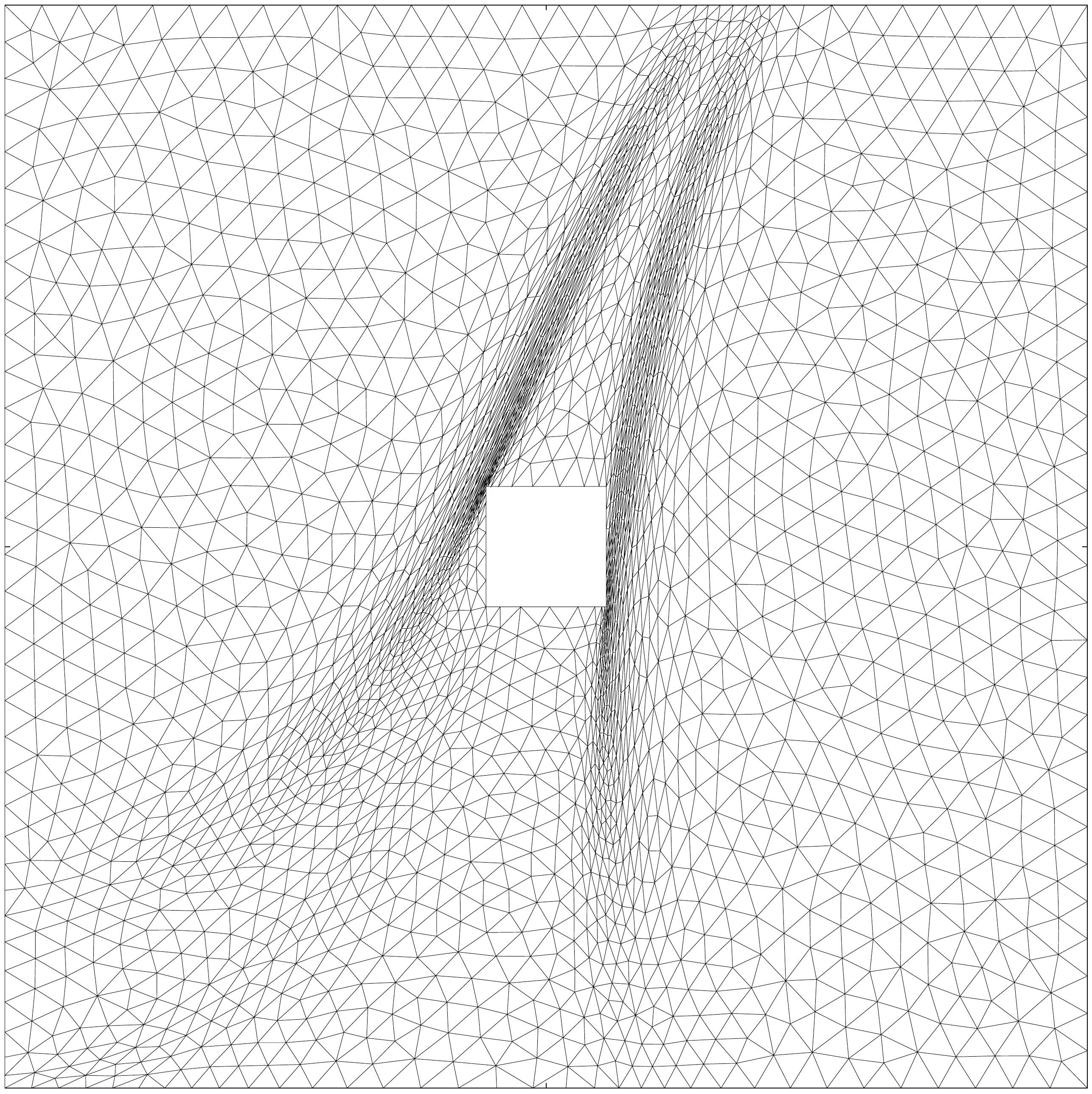}
      \includegraphics[width=0.23\textwidth,clip]{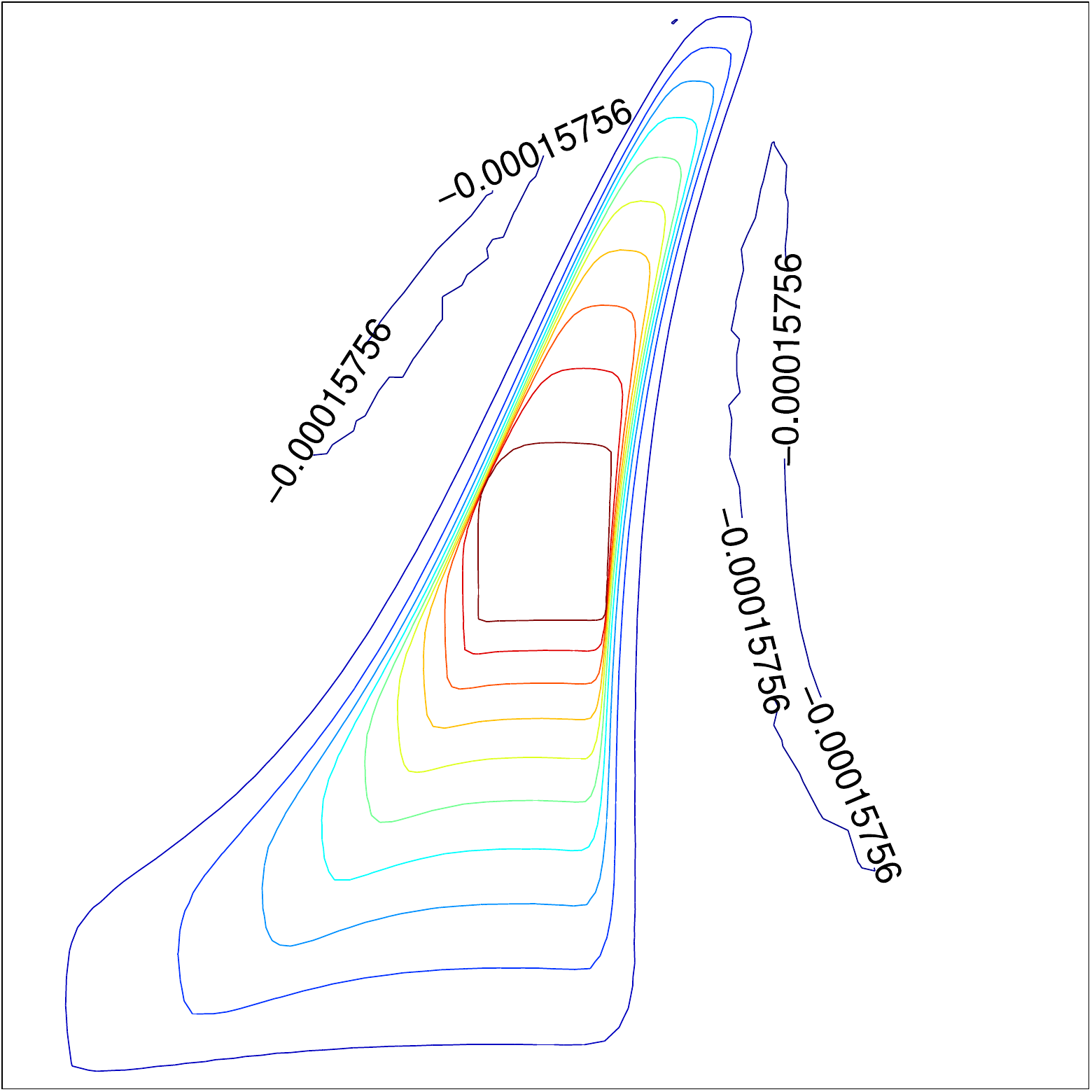}
      \label{fig:aniso1hb5}
   } 

   \caption{Anisotropic diffusion: meshes and contour plots of the numerical solution for the
      \protect\subref{fig:aniso1hb1} isotropic and
      \protect\subref{fig:aniso1hb5} HBEE-based ($M_{HB}$) metric tensors.
    }
      \label{fig:aniso1hb}
\end{figure}

\begin{figure}[t] \centering
   \subfloat[$M_{DMP}$: 4253 triangles, no undershoots, maximum aspect ratio 76.5.]{
      \includegraphics[width=0.23\textwidth,clip]{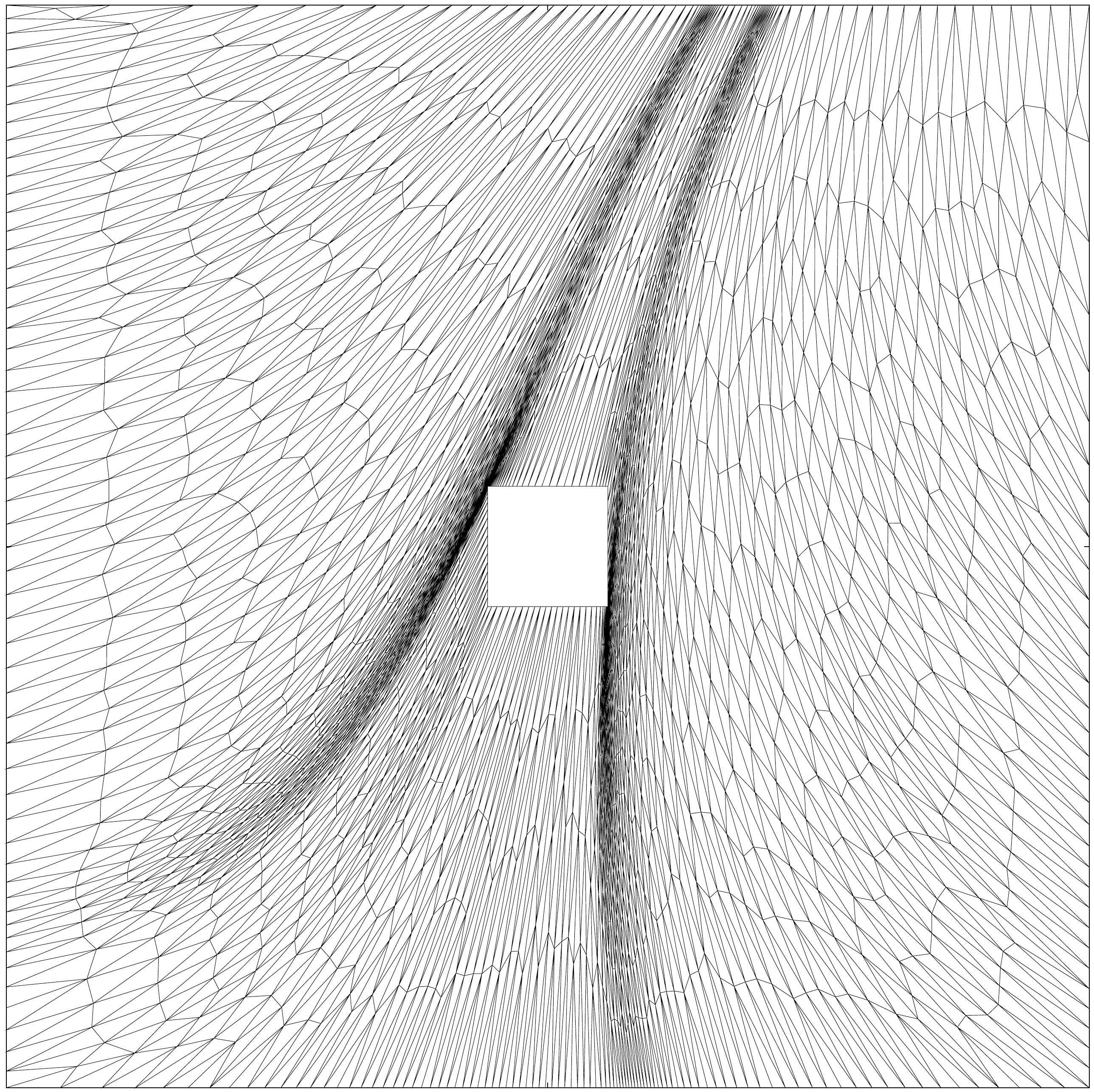}
      \includegraphics[width=0.23\textwidth,clip]{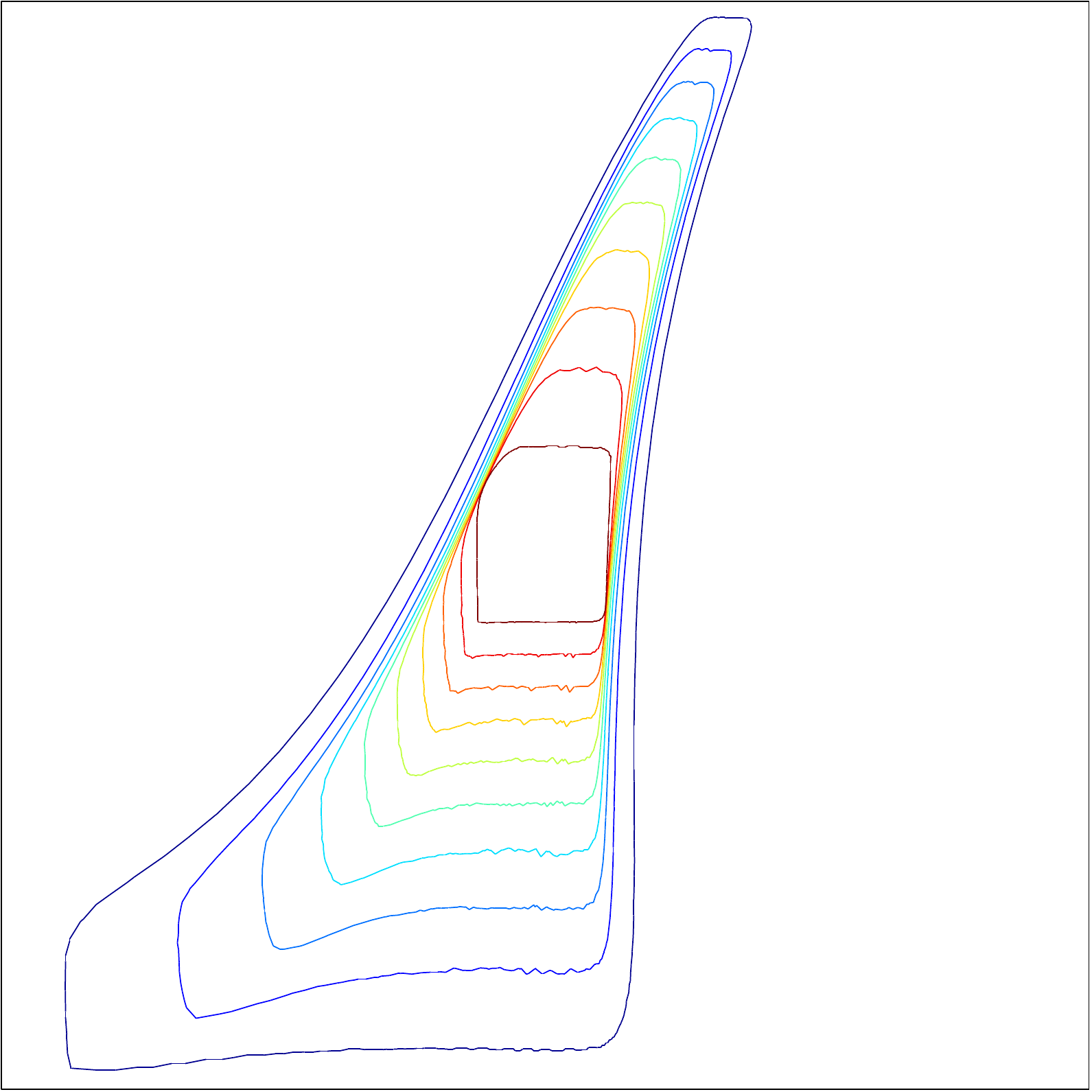}
      \label{fig:aniso1hb111}
   } \\
   \subfloat[$M_{DMP+HB}$: 4381 triangles, no undershoots, maximum aspect ratio 84.2.]{
      \includegraphics[width=0.23\textwidth,clip]{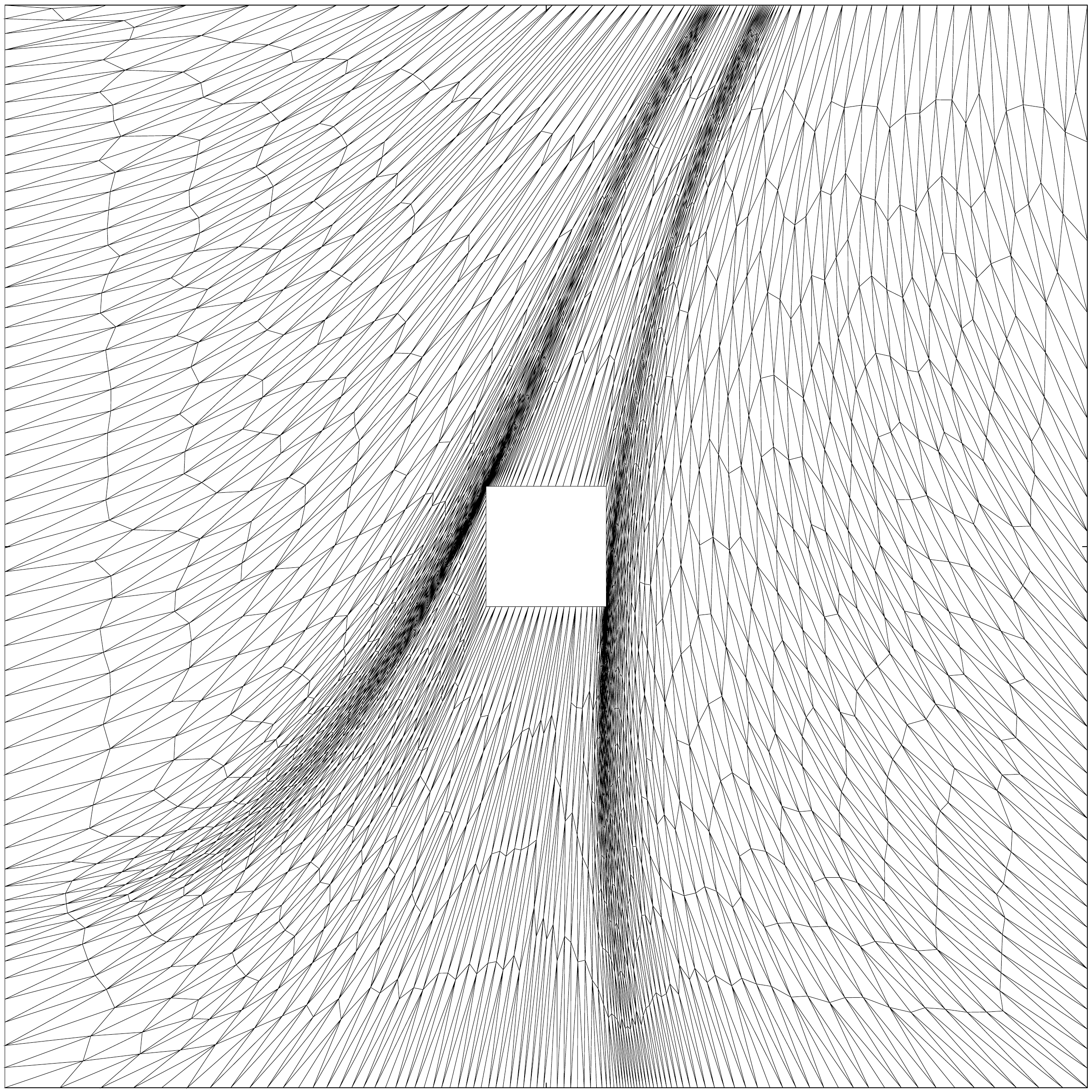}
      \includegraphics[width=0.23\textwidth,clip]{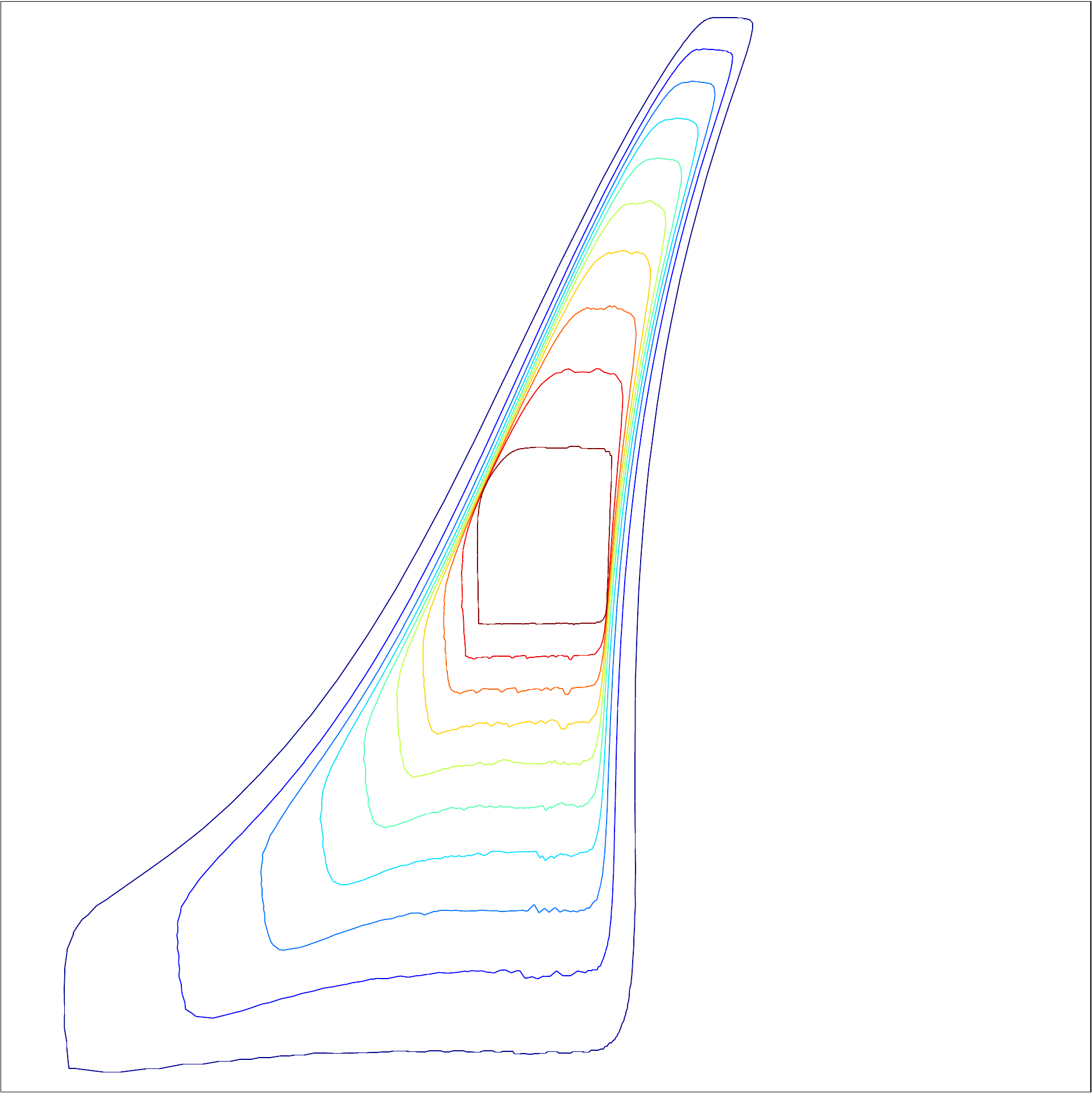}
      \label{fig:aniso1hb116}
   } 
  \caption{Anisotropic diffusion: meshes and contour plots of the numerical solution 
   for the DMP-compliant
  \protect\subref{fig:aniso1hb111} Hessian recovery-based ($M_{DMP+H}$) and
  \protect\subref{fig:aniso1hb116} HBEE-based ($M_{DMP+HB}$) metric tensors.
   }
  \label{fig:aniso1hb11}
\end{figure}

To emphasize the compliance with the DMP, the metric tensors $M_{DMP+H}$ and $M_{DMP+HB}$ based on \eqref{eq:B} and \eqref{eq:Bhbee}, respectively, are also compared to a uniform mesh and the anisotropic metric tensor $M_{HB}$ based solely on the HBEE (Sect.~\ref{sec:cornerSingularity}).
Figures~\ref{fig:aniso1hb} and \ref{fig:aniso1hb11} show meshes and solution contours. 

% anisotropic diffusion 

No overshoots in the finite element solutions are observed for all cases, but Fig.~\ref{fig:aniso1hb} shows that undershoots and unphysical minima occur in the solutions obtained with the uniform mesh ($\min u_h \approx -0.059$) and $M_{HB}$ ($\min u_h \approx -0.0032$).
As expected, no undershoots can be observed for $M_{DMP+H}$ and $M_{DMP+HB}$ (Fig.~\ref{fig:aniso1hb11}).

As for $M_{DMP}$ and $M_{DMP+HB}$, the solution contours for both metric tensors are almost undistinguishable and, although not as smooth, still quite comparable to the one obtained with $M_{HB}$ (cf. Figs.~\ref{fig:aniso1hb11} and \ref{fig:aniso1hb5}), thus providing a good adaptation to the sharp solution jump near the interior boundary (cf. the somewhat smeared solution jump for the isotropic mesh  in Fig.\ref{fig:aniso1hb1}, right). 
The mesh computed by means of HBEE is fully comparable with the mesh obtained using Hessian recovery. 
However, the maximum aspect ratio of the mesh obtained by means of the HBEE is, again, slightly larger.

\section{Concluding remarks}
\label{sec:concludingRemarks}
%********************************************************************
Numerical results show that a global HBEE can be a successful alternative to Hessian recovery in mesh adaptation: the new method performs successfully and is quite comparable with the commonly used Hessian recovery method.
A fast approximate solution is as fast as Hessian recovery and proved to be sufficient to provide enough directional information for the purpose of the mesh adaptation.
Moreover, as already observed in \cite{HuKaLa10}, it could be of advantage  for problems with discontinuities, because Hessian recovery could result in unnecessarily high mesh density for such problems.
The global HBEE seems to be less affected by this issue and, depending on the underlying problem, can provide a more appropriate mesh adaptation.
Also, Hessian recovery can also cause a light mesh smoothing: meshes obtained by means of Hessian recovery-based method in Sect.~\ref{sec:variationalProblem} have a smaller maximum aspect ratio than meshes obtained with HBEE and therefore seem to be slightly worse in terms of adaptation with the steep boundary layers.
For the corner singularity problem, the discussed mesh adaptation method proved to be able to capture the required information and to provide the optimal mesh grading without any \emph{a~priori} information on the solution.

One of the key components of the method is the reliability of the error estimator on anisotropic meshes: error estimation with hierarchical bases is usually based on the saturation assumption, which basically states that quadratic approximations provide finer information on the solution than linear ones.
Existing results on its validity require bounds on the elements' aspect ratio \cite{DoeNoc02}.
It is still unclear if similar results can be achieved for general adaptive meshes, but numerical results suggest that aspect ratio bounds are not necessary if the mesh is properly aligned.
Moreover, it seems that good mesh adaptation does not require an accurate Hessian recovery or an accurate error estimator, but rather some additional information of global nature, although it is still unclear which information exactly is necessary.
This  question is definitely of interest for future investigations.

%********************************************************************
\vskip 0.75em
\noindent
\textbf{Acknowledgements.} The author is very grateful to the anonymous referees for their valuable comments and suggestions.

%*** Bibliography ***************************************************

\small{

}

\begin{thebibliography}{10}
\providecommand{\url}[1]{{#1}}
\providecommand{\urlprefix}{URL }
\expandafter\ifx\csname urlstyle\endcsname\relax
  \providecommand{\doi}[1]{DOI~\discretionary{}{}{}#1}\else
  \providecommand{\doi}{DOI~\discretionary{}{}{}\begingroup
  \urlstyle{rm}\Url}\fi
\bibitem{AgLiVa08}
Agouzal, A., Lipnikov, K., Vassilevski, Y.: Generation of quasi-optimal meshes
  based on a posteriori error estimates.
\newblock In: Proceedings of the 16th International Meshing Roundtable, pp.
  139--148 (2008)

\bibitem{AgLiVa09}
Agouzal, A., Lipnikov, K., Vassilevski, Y.: Anisotropic mesh adaptation for
  solution of finite element problems using hierarchical edge-based error
  estimates.
\newblock In: Proceedings of the 18th International Meshing Roundtable, pp.
  595--610 (2009)

\bibitem{AgLiVa10}
Agouzal, A., Lipnikov, K., Vassilevski, Y.: Hessian-free metric-based mesh
  adaptation via geometry of interpolation error.
\newblock Comput. Math. Math. Phys. \textbf{50}(1), 124--138 (2010)

\bibitem{AgoVas10}
Agouzal, A., Vassilevski, Y.V.: Minimization of gradient errors of piecewise
  linear interpolation on simplicial meshes.
\newblock Comput. Methods Appl. Mech. Engrg. \textbf{199}(33-36), 2195--2203
  (2010)

\bibitem{Apel99}
Apel, T.: Anisotropic Finite Elements: Local Estimates and Applications.
\newblock B. G. Teubner, Stuttgart (1999)

\bibitem{ApGrJM04}
Apel, T., Grosman, S., Jimack, P.K., Meyer, A.: A new methodology for
  anisotropic mesh refinement based upon error gradients.
\newblock Appl. Numer. Math. \textbf{50}(3-4), 329--341 (2004)

\bibitem{BanSmi93}
Bank, R.E., Smith, R.K.: A posteriori error estimates based on hierarchical
  bases.
\newblock SIAM J. Numer. Anal. \textbf{30}(4), 921--935 (1993)

\bibitem{BanXu03}
Bank, R.E., Xu, J.: Asymptotically exact a posteriori error estimators, {P}art
  {I}: Grids with superconvergence.
\newblock SIAM J. Numer. Anal. \textbf{41}(6), 2294--2312 (2003)

\bibitem{BanXu03a}
Bank, R.E., Xu, J.: Asymptotically exact a posteriori error estimators, {P}art
  {II}: General unstructured grids.
\newblock SIAM J. Numer. Anal. \textbf{41}(6), 2313--2332 (2003)

\bibitem{Bildhauer03}
Bildhauer, M.: Convex variational problems. Linear, nearly linear and
  anisotropic growth conditions., \emph{Lecture Notes in Mathematics}, vol.
  1818.
\newblock Springer Berlin / Heidelberg (2003)

\bibitem{BoGHLS97}
Borouchaki, H., George, P.L., Hecht, F., Laug, P., Saltel, E.: Delaunay mesh
  generation governed by metric specifications. {P}art {I}. {A}lgorithms.
\newblock Finite Elements in Analysis and Design \textbf{25}(1-2), 61--83
  (1997)

\bibitem{BoGeMo97}
Borouchaki, H., George, P.L., Mohammadi, B.: Delaunay mesh generation governed
  by metric specifications. {P}art {II}. {A}pplications.
\newblock Finite Elem. Anal. Des. \textbf{25}(1-2), 85--109 (1997)

\bibitem{CaHuRu01}
Cao, W., Huang, W., Russell, R.D.: Comparison of two-dimensional r-adaptive
  finite element methods using various error indicators.
\newblock Math. Comput. Simulation \textbf{56}(2), 127--143 (2001)

\bibitem{CaHeMP97}
Castro-D{\'\i}az, M.J., Hecht, F., Mohammadi, B., Pironneau, O.: Anisotropic
  unstructured mesh adaption for flow simulations.
\newblock Int. J. Numer. Meth. Fluids \textbf{25}(4), 475--491 (1997)

\bibitem{CiaRav73}
Ciarlet, P.G., Raviart, P.A.: Maximum principle and uniform convergence for the
  finite element method.
\newblock Comput. Methods Appl. Mech. Engrg. \textbf{2}(1), 17--31 (1973)

\bibitem{DeLeYs89}
Deuflhard, P., Leinen, P., Yserentant, H.: Concepts of an adaptive hierarchical
  finite element code.
\newblock Impact Comput. Sci. Engrg. \textbf{1}(1), 3--35 (1989)

\bibitem{DoGrPf99}
Dobrowolski, M., Gr{\"a}f, S., Pflaum, C.: On a posteriori error estimators in
  the finite element method on anisotropic meshes.
\newblock Electron. Trans. Numer. Anal. \textbf{8}, 36--45 (1999)

\bibitem{Dolejs98}
Dolej\v{s}{\'\i}, V.: Anisotropic mesh adaptation for finite volume and finite
  element methods on triangular meshes.
\newblock Comput. Vis. Sci. \textbf{1}(3), 165--178 (1998)

\bibitem{DoeNoc02}
D{\"o}rfler, W., Nochetto, R.H.: Small data oscillation implies the saturation
  assumption.
\newblock Numer. Math. \textbf{91}, 1--12 (2002)

\bibitem{ForPer01}
Formaggia, L., Perotto, S.: New anisotropic a priori error estimates.
\newblock Numer. Math. \textbf{89}(4), 641--667 (2001)

\bibitem{Frey08}
Frey, P.J., George, P.L.: Mesh Generation. Second Edition.
\newblock John Wiley \& Sons, Inc., Hoboken, NJ (2008)

\bibitem{bamg}
Hecht, F.: BAMG: Bidimensional Anisotropic Mesh Generator (2006).
\newblock \\\url{http://www.ann.jussieu.fr/~hecht/ftp/bamg/}

\bibitem{Huang05a}
Huang, W.: Metric tensors for anisotropic mesh generation.
\newblock J. Comput. Phys. \textbf{204}(2), 633--665 (2005)

\bibitem{HuKaLa10}
Huang, W., Kamenski, L., Lang, J.: A new anisotropic mesh adaptation method
  based upon hierarchical a posteriori error estimates.
\newblock J. Comput. Phys. \textbf{229}(6), 2179--2198 (2010)

\bibitem{HuKaLi10}
Huang, W., Kamenski, L., Li, X.: Anisotropic mesh adaptation for variational
  problems using error estimation based on hierarchical bases.
\newblock Canad. Appl. Math. Quart. \textbf{17}(3), 501--522, arXiv:1006.0191
  (2009)

\bibitem{HuaLi10}
Huang, W., Li, X.: An anisotropic mesh adaptation method for the finite element
  solution of variational problems.
\newblock Finite Elem. Anal. Des. \textbf{46}(1-2), 61--73 (2010)

\bibitem{HS01}
Huang, W., Sun, W.W.: Variational mesh adaptation {II}: error estimates and
  monitor functions.
\newblock J. Comput. Phys. \textbf{184}(2), 619--648 (2003)

\bibitem{Johnson87}
Johnson, C.: Numerical Solution of Partial Differential Equations by the Finite
  Element Method.
\newblock Cambridge University Press (1987)

\bibitem{Kamens09PhD}
Kamenski, L.: Anisotropic mesh adaptation based on hessian recovery and a
  posteriori error estimates.
\newblock Ph.D. thesis, TU Darmstadt (2009)

\bibitem{LiHua10}
Li, X., Huang, W.: An anisotropic mesh adaptation method for the finite element
  solution of heterogeneous anisotropic diffusion problems.
\newblock J. Comput. Phys. \textbf{229}(21), 8072--8094 (2010)

\bibitem{Ovall07}
Ovall, J.S.: Function, gradient, and {H}essian recovery using quadratic
  edge-bump functions.
\newblock SIAM J. Numer. Anal. \textbf{45}(3), 1064--1080 (2007)

\bibitem{VaMDDG07}
Vallet, M.G., Manole, C.M., Dompierre, J., Dufour, S., Guibault, F.: Numerical
  comparison of some {H}essian recovery techniques.
\newblock Int. J. Numer. Methods Engrg. \textbf{72}(8), 987--1007 (2007)

\bibitem{VasLip99}
Vassilevski, Y., Lipnikov, K.: An adaptive algorithm for quasioptimal mesh
  generation.
\newblock Comput. Math. Math. Phys. \textbf{39}(9), 1468--1486 (1999)

\bibitem{ZhaNag05}
Zhang, Z., Naga, A.: A new finite element gradient recovery method:
  Superconvergence property.
\newblock SIAM J. Sci. Comput. \textbf{26}(4), 1192--1213 (2005)

\bibitem{ZieZhu92}
Zienkiewicz, O.C., Zhu, J.Z.: The superconvergent patch recovery and a
  posteriori error estimates. {P}art 1: The recovery technique.
\newblock Int. J. Numer. Methods Engrg. \textbf{33}(7), 1331--1364 (1992)

\bibitem{ZieZhu92a}
Zienkiewicz, O.C., Zhu, J.Z.: The superconvergent patch recovery and a
  posteriori error estimates. {P}art 2: Error estimates and adaptivity.
\newblock Int. J. Numer. Methods Engrg. \textbf{33}(7), 1365--1382 (1992)

\end{thebibliography}
\end{document}